# Benders, Nested Benders and Stochastic Programming: An Intuitive Introduction

James Murphy

# Benders, Nested Benders and Stochastic Programming

An Intuitive Introduction

**James Murphy**

December 2013

## Abstract

This article aims to explain the Nested Benders algorithm for the solution of large-scale stochastic programming problems in a way that is intelligible to someone coming to it for the first time. In doing so it gives an explanation of Benders decomposition and of its application to two-stage stochastic programming problems (also known in this context as the L-shaped method), then extends this to multi-stage problems as the Nested Benders algorithm. The article is aimed at readers with some knowledge of linear and possibly stochastic programming but aims to develop most concepts from simple principles in an understandable way. The focus is on intuitive understanding rather than rigorous proofs.

## Keywords

Benders, Nested Benders, stochastic programming, linear programming, constrained optimization, L-shaped method, decomposition methods, cutting-plane methods.

# Contents

# Introduction

The aim of this article is to give an explanation of the *Nested Benders* algorithm that is intelligible to someone approaching it for the first time, equipped with some basic knowledge of linear programming and possibly stochastic programming. The article is a tutorial, rather than a literature review. As such, the references have been kept fairly limited and no claim is made that those papers and books cited are the latest, best, most complete, etc.

The article will develop the Nested Benders algorithm for stochastic programming problems in the following way. First, two asides, the first on linear programming and the second on Duality theory, will present some of the basic ideas from these two topics that are necessary in the development of the Benders decomposition algorithm. Readers familiar with these subjects can skip these asides, referring to them as necessary. Following this, the basic ideas of Benders decomposition for constrained optimization problems will be introduced. Up to this point the article deals with the Benders decomposition method in general without specific reference to stochastic programming. The subsequent parts of the article apply Benders decomposition to stochastic programming problems. These will be examined in more detail to show how their special structure lends itself to solution by decomposition methods. The details of how decomposition leads to an efficient algorithm for the solution of two-stage stochastic programming problems will then be explained. This algorithm is then extended to multi-stage stochastic programming problems. At this point, we will have arrived at our goal, the Nested Benders algorithm for stochastic programming problems. Finally the article will conclude with a very brief discussion of some additional topics such as parallelization and a list of alternative methods.

Nested Benders is an important algorithm in the solution of multi-stage stochastic programming problems, which, by their nature, tend to be large, necessitating the use of efficient solution algorithms. It is based on the idea that problems can be decomposed into their component parts, which can be tackled somewhat independently and then reintegrated in order to solve the original, more complex problem. Stochastic programming problems lend themselves to this type of solution because they can be thought of as trees, with each node in the tree representing a decision point. At each one of these decision points an optimization problem needs to be solved and the solution of all these sub-problems then needs to be integrated to come up with an optimal decision at the root node. Of course, the decisions taken in one time step affect the problems that need to be solved in the next one. Such problems are often known as *recourse problems* because an initial decision is made, time elapses and events occur and then a further decision is made (a recourse decision) that aims to be optimal in light of the events that have just occurred. We will return to this and to stochastic programming in a later section. First, however, it will be helpful to develop the basic ideas of Benders decomposition without the extra confusion and terminology introduced by stochastic programming.

Before starting to explain the Benders algorithm, some basic ideas from linear programming and duality theory should be introduced. Readers familiar with this material can safely skip it, referring to it if necessary.

## Aside on Linear Programming

This aside on linear programming cannot possibly be a complete introduction to the subject. It merely aims to introduce some of the ideas about solving linear programmes that are essential for the development of the Benders algorithm. These are the ideas of *feasibility* and *boundedness* of linear programmes and the fact that an optimal solution to a linear programme always lies at a corner of its feasible area. We will not prove this last point, only make an argument for its reasonableness. For further information on the very large subject and theory of linear programming there are many good introductory texts available (e.g. Vanderbei (2000), my personal favourite), and the interested reader should consult one or more of them.

### What is Linear Programming?

Firstly, if you are unsure about what linear programming is then this article is probably not the right place for you to start and you should consult an introductory text as mentioned above. Only a very brief introduction is given here in order to standardize terminology and to give a reference for the rest of the article.

Linear programming is the subject of finding optimal values of linear functions of a number of variables subject to linear constraints on those variables. It is a form of constrained optimization that is particularly interesting because it is tractable (i.e. it is actually possible to find solutions to such problems in a reasonable amount of time) and a surprisingly large number of problems can be cast (or approximated) as linear programming problems.

Stated mathematically linear programmes have the form

$$\text{maximize } \sum_{j=1}^{n} c_j x_j$$
$$\text{subject to } \sum_{j=1}^{n} a_{ij} x_j \leq b_i \qquad i = 1,\ldots,m \qquad (1)$$
$$x_j \geq 0 \qquad j = 1,\ldots,n$$

where here there are *n* variables to be optimized (the $x_j$'s) and *m* constraints on those variables. The objective function (top line) is a linear combination of the variables, as are the constraints (second line). Note that this can easily be converted to a minimization problem by negating the objective function. The problem (1) can be re-written in matrix-vector form as

$$\begin{aligned} \max_{x} \quad & c'x \\ \text{s.t.} \quad & Ax \leq b \\ & x \geq 0 \end{aligned} \qquad (2)$$

where here $x$ and $c$ are $n \times 1$ vectors, $b$ is a $m \times 1$ vector and $A$ is a $m \times n$ matrix. The abbreviation "s.t." is very commonly used for the term "subject to" (or "such that").

It might appear that (1) is not really a general form for a linear programme, as we cannot deal with any variable that takes values less than 0. This appears to be a severe limitation as there are many problems where we might want variables to be able to take negative values or even to take any value from the entire set of real numbers, both positive and negative. Happily, it is possible to easily fit such variables into problems with the form of (1) with just a little bit of pre-processing.

In order to deal with a non-positive variable, we simply change the sign in front of every occurrence of that variable in the constraints and objective. We can then insist that it be positive. The only thing to remember is to take the negation of the optimal value found for this variable when reporting the results of the original problem.

To deal with a variable that can take any real value we must split it up into positive and negative parts, such that both are non-negative and only one is non-zero at any time. The original variable $x$ is then replaced in all the constraints and the objective with the expression $(x^+ - x^-)$, formed of the two non-negative variables $x^+$ and $x^-$. When reporting the optimal value of x, we simply report the value of $(x^+ - x^-)$ for the optimal values of $x^+$ and $x^-$.

**Feasibility**

A linear programme is feasible if there is some value that the variables can take so that all the constraints are simultaneously satisfied. The set of values for which all the constraints are satisfied are called the feasible set (or feasible region).

For a very simple problem with two variables we can show this diagrammatically. The variables (call them $x_1$ and $x_2$) form the axes of the plot. The linear constraints can be shown as lines on the graph, values to one side of which are not permitted (each constraint defines a permitted and non-permitted half-plane). The area on the permitted side of all these constraints (or the intersection of all the permitted half-planes) is the feasible area.

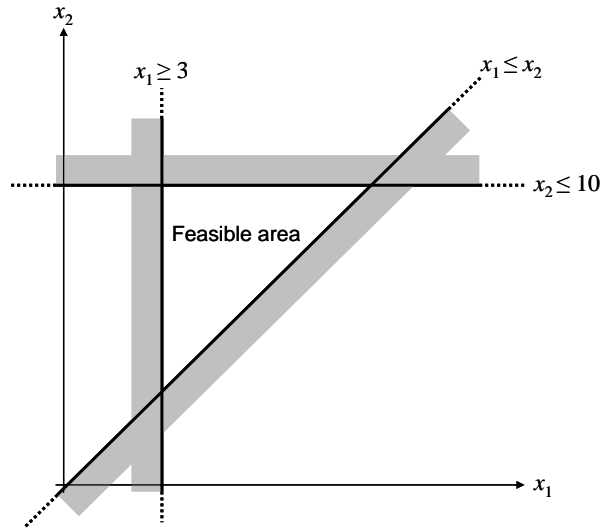

**Figure 1**  *Simple problem in two variables with three constraints.  Constraints are shown as lines with the grey shading indicating the side on which the variables are not permitted to lie.*

### Boundedness

A linear programme is bounded if it is feasible and if its (linear) objective is bounded within the feasible region (i.e. if it cannot be made to go to +∞ for a maximization or –∞ for a minimization).  This will always be the case if the feasible area is finite.  Even if the feasible area is semi-infinite the problem can still be bounded as the constraints might place an upper bound on the objective value (for a maximization, or a lower bound for a minimization).  The following figure illustrates this point.  Note that for a linear objective the objective function has a single direction in which it increases most quickly (i.e. it has a constant gradient everywhere).

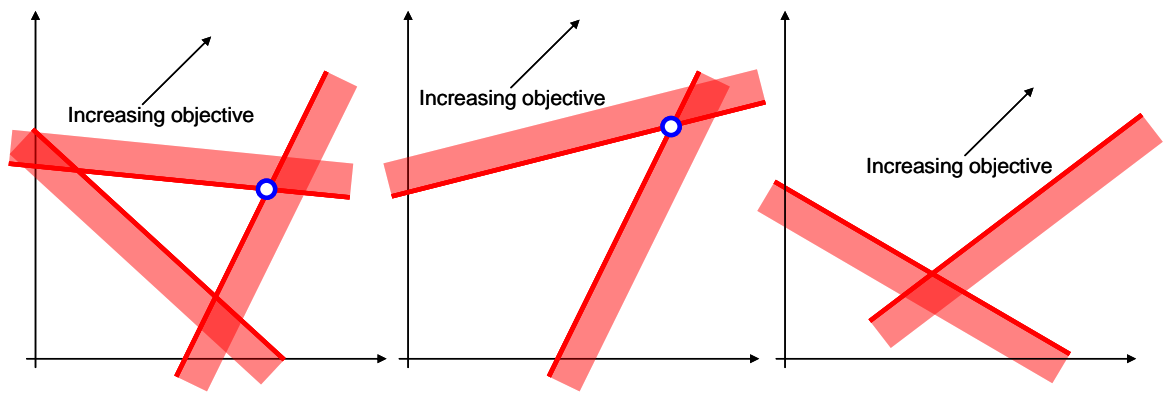

**Figure 2**  *Three examples of maximization.  The red lines are constraints and the objective increases quickest along x=y (e.g. x+y).  The blue circle shows the optimal point.  The leftmost example has a finite feasible region and is bounded.  The middle example has a semi-infinite feasible region but is still bounded.  The right example has semi-infinite feasible region and is unbounded.*

## Solution at Cornerpoint

In linear programming problems an optimal solution (if it exists, i.e. if the problem is feasible and bounded) will always lie at a cornerpoint of the feasible area. To see this we observe that the objective function is linear, meaning that it has a gradient that is constant everywhere. It increases fastest along this gradient. So, at any given point to increase the value of the objective function we will move up the gradient. Starting with a feasible point we can increase the objective function in such a way until we reach one of the constraints. If this constraint is perpendicular to the gradient then the objective function will be constant all the way along it (since, by definition, in the direction perpendicular to its gradient there is no change in the value of a function). In this case we have reached the optimal objective value and it lies all the way along this constraint (including at the points where this constraint meets other constraints, corners of the feasible area). This is true because this constraint forces the solution to not lie beyond the particular contour of the objective function that it defines (since it is perpendicular to the gradient it must define a contour of the objective function), and so the objective function cannot take a value greater than the one on the line of this constraint.

If the constraint is not perpendicular to the gradient then there will be a direction along it that we can move where the objective function will increase. Since the objective function is linear this rate of increase will be constant all the way along the length of the constraint. If we continue along the constraint (and the problem is bounded) we will eventually meet another constraint. This too will have a direction along it where the objective value will increase. Because of the preceding constraint, however, we can only follow this new constraint in one direction and remain in the feasible area. If this direction is the same as the one in which the objective value increases we can follow this new constraint to feasible points with higher objective values, eventually encountering another constraint and repeating the process. If, however, the direction in which we can follow the new constraint is the one in which the objective function decreases we should not go any further as we have reached the optimal objective value as there is no further direction to follow that will increase the objective. This is at the corner between the two constraints. So, in all cases an optimal point lies at a cornerpoint of the feasible area. In cases where the constraint on which it lies is not perpendicular to the objective gradient, the optimal point is uniquely at the corner.

The idea of following constraints around the edge of the feasible area until an optimal cornerpoint is reached is the basis of the simplex algorithm for solving linear programming problems. This algorithm (invented by George Dantzig in 1947 and first published in Dantzig (1951)) is crucially important in linear programming as it was the first efficient way of finding solutions to linear programming problems and is still widely used today.

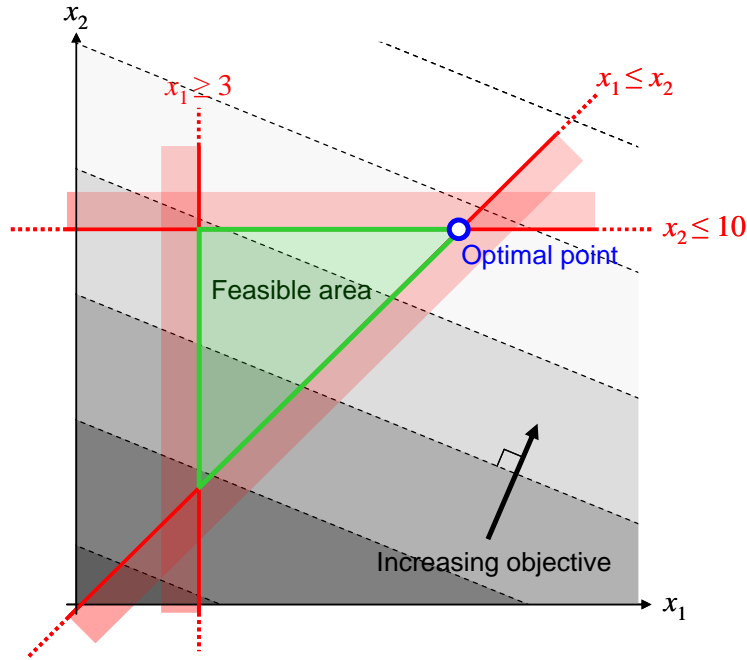

**Figure 3** *Illustration of an optimal solution always lying on a cornerpoint of the feasible area. Constraints are shown red. The grey shading indicates the increasing objective function with dotted lines being evenly spaced contours. No matter what further constraints were introduced the optimal point would always lie at a cornerpoint. The problem shown is a maximization.*

## Aside on Duality Theory

A key concept in linear programming is the idea of the *dual* of a linear programming problem. This is another linear programming problem that can be used to place upper bounds (for a maximization) on the value of the objective function. Every linear programme has a dual and the dual of the dual of a problem is the original (or *primal*) problem. The dual can also be used to discover information about the boundedness and feasibility of the primal problem. This section will briefly introduce those parts of duality theory that we need to develop the methods described in the rest of this article. It is a necessarily brief treatment of the subject and no results are proved. For a fuller treatment consult e.g. Vanderbei (2000).

### Formation of Duals

Any feasible solution to a maximization problem (we will deal exclusively with maximization henceforth in this section) obviously gives a lower bound on the optimal objective value. The dual of a linear programming problem can be seen as arising out of the search for an upper bound on the optimal objective value.

Consider the maximization problem in (1). The objective function is composed of a linear combination of the problem variables, the $x$'s. How can we place an upper bound on this? The only parts of the problem that gives us an upper bound on anything are the constraints in the second line. These give an upper bound on various linear combinations of the variables. Since all of the *x*'s are positive, we can use these to give us an upper bound on the objective function.

If we take a weighted sum (with non-negative weights) of these constraints in such a way that the coefficients of each of the variables are at least as large as those in the objective function we will, since the variables are positive, have placed a bound on the objective function. If we call the multiplier of each constraint $y_i$ this weighted sum will give us

$$\sum_{i=1}^{m} y_i \sum_{j=1}^{n} a_{ij} x_j \leq \sum_{i=1}^{m} y_i b_i \qquad (3)$$

along with the non-negativity condition on the $y_i$'s

$$y_i \geq 0 \qquad\qquad i = 1,\ldots,m. \qquad (4)$$

From the left hand side of (3) we can get an expression for the coefficient of each variable $x_j$, which we require to be at least as large as the coefficient of that $x_j$ in the objective function of the primal problem $c_j$. This gives the following set of *n* constraints, one for each of the original variables.

$$\sum_{i=1}^{m} y_i a_{ij} \geq c_{ij} \qquad\qquad j = 1,\ldots,n \qquad (5)$$

To get the tightest upper bound we can, we want to choose the values of $y_i$ in such a way as to minimize the right hand side of (3), since this is our bound on the optimal objective. This along with the conditions in (4) and (5) gives us the following problem.

$$\begin{aligned} \text{minimize} \quad & \sum_{i=1}^{m} y_i b_i \\ \text{subject to} \quad & \sum_{i=1}^{m} y_i a_{ij} \geq c_{ij} \qquad j = 1,\ldots,n \\ & y_i \geq 0 \qquad\qquad\quad i = 1,\ldots,m \end{aligned} \qquad (6)$$

This is another linear programming problem, this time with *m* variables and *n* constraints. It is the dual problem of (1) and since that was a general linear programming problem this shows us how to form the dual of any linear programming problem.

## Optimality

The optimal solution to a linear programme is bounded below by a feasible solution of the primal problem and above by a feasible solution of the dual problem. As we saw in the previous section, any feasible solution of the dual is always at least as great as any feasible solution of the primal (this is the *Weak Duality Theorem*).

The *Strong Duality Theorem* asserts that if the primal problem has any feasible solutions then the optimal value of the primal and dual problems is the same. This is the same as saying that there is one point of intersection between objective values for feasible dual solutions and objective values for feasible primal solutions. Clearly this point of intersection is the optimal objective value for both the primal and the dual problems, since it must be the maximum obtainable objective value for the primal problem (a maximization) and the minimum attainable objective value for the dual problem (a minimization).

Furthermore, this means that if the primal problem is feasible and bounded then the dual must also be feasible and bounded since it has at least one feasible point (the one that produces the optimal objective value) and its objective is bounded (below) by this.

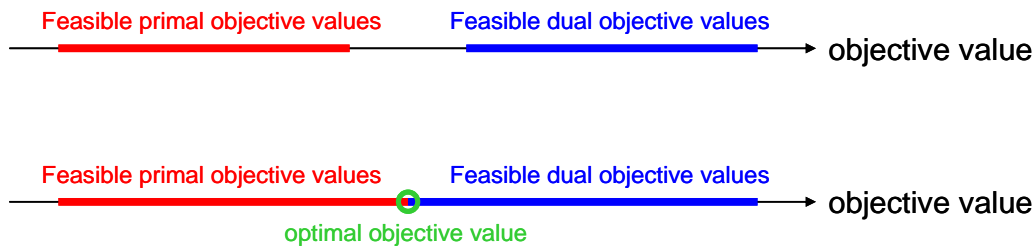

**Figure 4** *Illustration of weak duality theorem (top) and strong duality theorem (bottom). (Diagram idea from Vanderbei (2000)).*

The strong duality theorem will not be proved here, but can be used as a test for optimality in cases where we have a solution for both the primal and the dual problem (this can be the case in the simplex algorithm, for example).

### Relationships between Primal and Dual Feasibility and Boundedness

The weak duality theorem, which states that all feasible primal points give an objective value lower than the objective value obtained from any feasible dual point, allows us to use the dual problem to give us information about the primal problem's feasibility and boundedness and *vice versa*.

If the primal problem is unbounded then its feasible objective values stretch off to $+\infty$. Since the feasible dual objective values must all lie above all feasible primal objective values, there is nowhere that they can be and so the dual must be infeasible. Similarly, if the dual is unbounded the primal problem is infeasible. This gives us the following two implications:

- Primal unbounded $\Rightarrow$ dual infeasible,
- Dual unbounded $\Rightarrow$ primal infeasible.

We might hope that these hold in the other direction, but it is not always the case. It can be the case that both the primal and dual problems are infeasible (we will not show that here but Vanderbei (2000) gives a simple example, which is sufficient to prove that it is possible). So, if we have infeasibility in the dual, we know that the primal is either unbounded or itself infeasible (and similarly for the primal). This gives us

- Primal infeasible $\Rightarrow$ dual infeasible or unbounded,
- Dual infeasible $\Rightarrow$ primal infeasible or unbounded.

Finally, from the strong duality theorem we know that if the primal is feasible and bounded then the dual is also feasible and bounded and their optimal value is the same.

- Primal feasible and bounded $\Leftrightarrow$ dual feasible and bounded.

## Benders: Basic Idea

We have now covered enough preliminaries to start to introduce the main feature of this article, the Nested Benders algorithm. Nested Benders, as we shall see later, is based on an iterative application of *Benders decomposition* and so it is this that we first need to understand if we are going to reach a thorough understanding of Nested Benders. Decomposition in this context means the breaking up of a *mathematical programming problem* (also known as a *constrained optimization problem*, of which linear programming is a particularly interesting subset) into separate smaller problems that are hopefully easier to solve separately and then re-integrate to get an overall solution. This section on the basic idea of Benders decomposition is inspired by the introductory article "Benders' Method" by Jensen (2003).

Benders (1962) originally applied the idea of decomposition to mixed integer programming problems. These problems are constrained optimization problems in which some variables can take real values and others may only take integer values. Decomposition was used by Benders to split the problem into a pure integer programming problem and a pure linear programming problem that could each be solved iteratively in turn to arrive at the solution to the overall problem. The advantage of this scheme is that the two sub-problems are each easier to solve than the original problem and so, although they might have to be solved several times to arrive at the solution, this is still likely to be quicker than trying to solve the original large problem.

When we applying decomposition to stochastic programming problems, we will not be trying to split an integer programming problem from a linear programming one, but instead trying to split a very large linear programming problem (the *deterministic equivalent* of the stochastic programming problem – more on this later) into smaller linear programming problems that are each much quicker to solve than the original. However, in this section the exact form of the sub-problems does not matter too much. For stochastic programming we will be interested in the case when all the problems and sub-problems are linear programming problems.

We start with a problem of the following form, which we call P.

P:
$$[Z^* =] \quad \max c_1 x + c_2 y$$
$$\text{s.t.} \quad A_1 x + A_2 y \leq b$$
$$x \geq 0$$
$$y \geq 0$$
(7)

Here $x$ and $y$ are vectors of variables that we are trying to optimize over. The constraints on these variables are given by the constraint matrices $A_1$ and $A_2$, along with the vector $b$. Any linear programming problem like (1) could be written in this way, splitting the problem into two sub-problems, with different variables in each. In the original case used by Benders, the $y$ variables were integer variables and the $x$ variables could take real values.

Now, if we fix the value of $y$ (i.e. fix the values of every element variable of $y$) we obtain the following linear programming problem, which we call PX.

PX:
$$[v_0(y) =] \quad \max c_1 x + c_2 y$$
$$\text{s.t.} \quad A_1 x \leq b - A_2 y$$
$$x \geq 0$$
(8)

Here the $c_2 y$ term in the objective and the $A_2 y$ term on the right hand side of the constraints are constant terms, since the value of $y$ is now fixed. $v_0(y)$ is the optimal value of the objective function when the value of $y$ is fixed to $y$. The true maximum value of the objective of (7), $Z^*$ must be at least as large as this, with equality only holding when the value of $y$ is chosen to be optimal for the original problem ($y^*$).

$$Z^* \geq v_0(y)$$
$$Z^* = v_0(y^*)$$
(9)

So now we have expressed a lower bound on the optimal objective value for the problem $Z^*$. We can also take the dual of problem PX and, as we saw in the *Aside on Duality Theory*, this allows us to place an upper bound on the value of $Z^*$. This is useful since our original problem P aims to maximize $Z^*$. Placing a lower bound on its value does not really limit this optimization, but an upper bound provides a potentially useful limit on the solution of the original problem. Furthermore, this upper bound can be expressed as a function of $y$ and a new dual variable $u$, eliminating the dependence of this bound on $x$, which will prove useful later. In this way, information about solving an optimization problem with a fixed value of $y$ can be used to feed back information to the original optimization problem. The key step in being able to do this is taking the dual of problem PX, which we know has the same optimal objective value as PX (see the *Aside on Duality Theory*).

To simplify the taking of the dual of PX, we can remove the constant $c_2 y$ term from the objective. Since this is just a constant term (as $y$ is held constant in PX), this does not affect the optimization itself and simply reduces the optimal objective value by $c_2 y$.

The dual of PX (with the $c_2 y$ term removed from its objective) is given by

$$\text{D1:} \quad \begin{aligned} [w_0(y) =] \quad & \min \; u(b - A_2 y) \\ \text{s.t.} \quad & u A_1 \geq c_1 \\ & u \geq 0 \end{aligned} \tag{10}$$

As the optimal objective value is the same for a primal problem and its dual, the optimal objective value for D1, $w_0(y)$, must be equal to that for PX with $c_2 y$ removed from its objective. Therefore,

$$w_0(y) = v_0(y) - c_2 y. \tag{11}$$

As was mentioned in the *Aside on Linear Programming*, an optimal solution to a linear programming problem lies at a cornerpoint of its feasible region (if such a region exists and if the problem is bounded). Call the set of cornerpoints of D1's feasible region $U$ and note that this is a finite (and possibly empty) set. If we assume that D1 is bounded (which may not be the case if the problem PX is infeasible for the given value of $y$ – more on this below) then, since the optimal solution must lie at a cornerpoint of the feasible region,

$$w_0(y) = \min_{u_j \in U} u_j (b - A_2 y). \tag{12}$$

That is, the minimum value of the objective of (10) $w_0(y)$ is equal to the minimum value of the objective of (10) attained on any of the cornerpoints of that problem's feasible region (the $u_j$ in the set of all cornerpoints $U$). Substituting this into our expression (11) relating the optimal objectives of PX and D1, we get

$$v_0(y) = c_2 y + \min_{u_j \in U} u_j (b - A_2 y) \tag{13}$$

As we saw above in (9), optimizing $v_0(y)$ over $y$ will yield $Z^*$, the optimal solution of the original problem, and so such an optimization is equivalent to the original problem. Call this problem P*.

$$\text{P*:} \quad [Z^* =] \quad \max_y \left[ c_2 y + \min_{u_j \in U} u_j (b - A_2 y) \right] \tag{14}$$

Since P* involves maximizing over a finite minimization we can use a standard trick in linear programming to remove it from the objective. If we consider any value that $u_j$ could take we know that the value of the term in square brackets must be at most the value given by this chosen value of $u_j$, since we are looking for the minimum value of the second term over all values of $u_j$. We can thus place an upper bound on the optimal objective value (as a function of $y$) for each possible value of $u_j$. This allows us to remove the minimization from the objective term and form the following equivalent problem, $P^k$.

$$P^k: \quad \begin{aligned} &[Z^* =] \quad &&\max Z \\ &\text{s.t.} \quad &&Z \le c_2 y + u_1(b - A_2 y) \\ & &&Z \le c_2 y + u_2(b - A_2 y) \\ & &&\ldots \\ & &&Z \le c_2 y + u_k(b - A_2 y) \end{aligned} \quad (15)$$

Here we assume that there are $k$ cornerpoints in $U$ (i.e. $|U| = k$), so the minimization term in the objective of P* yields $k$ independent constraints on the value of $Z$. These upper bounds are functions of $y$ (and the $u_j$'s); the $x$ variables from the sub-problem have been completely eliminated. Since the $u_j$'s can be dealt with by enumeration, the problem $P^k$ (which as we have seen is equivalent to the original problem P) can be dealt with purely in terms of the $y$ variables. In the original application of Benders' algorithm to mixed integer programming problems this meant that the original problem was reduced to a pure integer programming problem in $y$.

So far, we are not really any better off than we were when we first considered P. Although we have eliminated $x$ from $P^k$ it could take considerable effort to compute the $u_j$'s and it is not clear that $P^k$ will be particularly easy to solve.

However, we can form a problem from $P^k$ that should be easier to solve by dropping some of the constraints from the problem (this is known as a relaxation of the problem). Consider the problem formed by only including the first $r$ constraints of the $k$ that come from the minimization over the cornerpoints of $U$. This problem is sometimes known as the *master problem*, to which constraints are added. Call this problem $P^r$.

$$P^r: \quad \begin{aligned} &[Z^r =] \quad &&\max Z \\ &\text{s.t.} \quad &&Z \le c_2 y + u_1(b - A_2 y) \\ & &&Z \le c_2 y + u_2(b - A_2 y) \\ & &&\ldots \\ & &&Z \le c_2 y + u_r(b - A_2 y) \end{aligned} \quad (16)$$

Now, if we solve this relaxed problem it may be that the solution we find violates one of the $k-r$ constraints that we have removed, in which case it is not a solution to P and we are no better off than before. On the other hand, the solution to the relaxed problem could be the solution to P and then we have succeeded in solving P by solving the easier relaxed problem. The trouble is, how can we tell?

Since $Z^r$ is the optimal solution to a less constrained maximization problem than $P^k$, the optimal solution to $P^k$ (and hence to P) must be bounded above by $P^k$ (since if there was a better optimal value that was feasible for the more constrained $P^k$ it would have been found as an optimal solution to $P^r$). So, we have

$$Z^* \le Z^r \qquad\qquad r \le k. \quad (17)$$

And, from (9) we also have a lower bound on the optimal value of Z*. Combining this with (11) allows us to express this lower bound in terms of the optimal solution to the dual of the sub-problem:

$$v_0(y) = w_0(y) + c_2 y \leq Z^* \tag{18}$$

If the upper and lower bounds for Z* are equal then we have found the optimal value of Z*. This leads us to our first version of the Benders algorithm.

### Algorithm

0. (Initialization) Choose an initial value for *y*, set *r* = 1.
1. Solve D1 using the current value of *y*.
2. Solve P$^r$ using the current value of *r*.
3. Compare $Z^r$ (obtained in step 3) with $w_0(y) + c_2 y$ (obtained from step 2).
    a. If they are equal then stop. $Z^* = Z^r$ and $y^r$ gives the optimal values for *y*.
    b. If they are not equal then set $y = y^r$ and increment *r* by 1 (i.e. add an extra constraint to the problem P$^r$). Return to step 1.

Since P$^k$ is equivalent to P, the algorithm will terminate at the latest when *r* = *k*, since at this stage step 2 will be solving the complete optimization problem. The hope is that the algorithm will in fact terminate much earlier than this. Early termination is achieved by judicious choice of the constraint to add to the problem P$^r$ in step 3a. Coming up with a way to sensibly add these constraints for maximum effect will be the focus of the section on *Applying Benders to Stochastic Programming Problems*. However, this section should have given some insight into the basic idea of problem decomposition and how large optimization problems can be solved by splitting them up into sub-problems. We will see in the next section that stochastic programming problems lend themselves to such decomposition and so make ideal problems on which to apply Benders decomposition and variants of it.

### Feasibility and Boundedness

In the above we assumed that the dual problem D1 was both feasible and bounded. This might not be the case. For clarity, this was ignored above, but it must now be addressed.

In fact, we can safely ignore the case where D1 is infeasible because D1 is the dual of PX. If D1 is infeasible then from duality theory (see *Aside on Duality Theory*) we know that the primal problem PX is either unbounded or itself has no feasible solution. Since the area of feasibility of D1 does not depend on the value of *y*, the infeasible or unbounded nature of PX applies for all values of *y*. This means that it is not possible to choose a value of *y* to make PX well-behaved and so the original problem P is also either infeasible or unbounded as solving it requires solving PX for some (optimal) value of *y*. Thus, the original problem can be condemned as uninteresting and we can safely ignore it (though, of course, any algorithm must be able to detect it).

If D1 is unbounded then it means that PX is infeasible for the given value of *y* (this is highly possible as, for example, the value of *y* is initially chosen arbitrarily). This is useful information, but in order to complete step 3 of the algorithm above we need to come up with a solution for D1. If it is unbounded, we can do this by adding in a constraint that will provide a bound for the solution in all circumstances. For example, the following would suffice, where *M* is a large positive number and *e* is a vector of 1's of the same dimensionality as *u*.

$$e'u \leq M \tag{19}$$

In the sections below we will solve this problem more satisfactorily, generating a useful constraint from the information about infeasibility. Here it is sufficient to realise that such an arbitrary additional constraint is justified (if not optimal) because once we know that the problem PX is infeasible for the current value of *y*, we know that this value of *y* will be rejected and a new one will be found after the addition of more constraints to the P$^r$ problem.

## Stochastic Programming Problem Structure

The key question raised by the algorithm in the previous section is how to add constraints to the problem in such a way as to improve the chance of arriving at the optimum solution before all the constraints are added.

At this point we can start to think about the structure of the problems we are trying to solve and whether this will help us to formulate our methods of dividing problems into sub-problems and then of adding constraints to the master problem in some sensible way.

In this article we will be interested in solving *stochastic linear programming* problems. These are *multi-stage* linear programmes that take account of uncertainty in at least some of the quantities involved in the problem. By multi-stage we mean problems in which an optimal initial decision is made, more information becomes available and then further decisions are made. (The term *decision* here is used to mean the solution to an optimization problem, specifically the optimal values of the variables of the problem, but this will be used interchangeably with the term *solution*). If the problems at each stage are linear programming problems then the problem is a multi-stage linear programme. If there is some uncertainty about the way in which the information can evolve from one period to the next and this is taken into account in decision making, the problem is known as a *stochastic* programming problem.

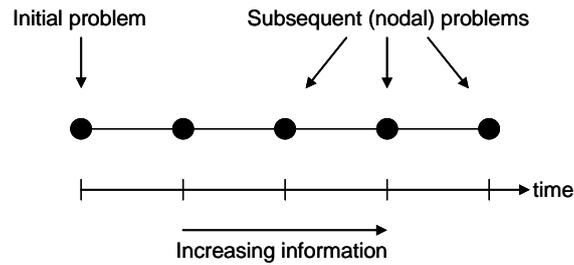

**Figure 5** *Example of a multi-stage problem with five stages. Each black dot represents a decision problem solved in light of the information available at that time and the decisions made at previous times. There is no uncertainty about the evolution of information, so this is a* deterministic dynamic *problem.*

We will see that stochastic programming problems lend themselves naturally to decomposition and understanding this natural decomposition will help us to come up with intuitive ways of adding constraints. Readers familiar with the idea of deterministic equivalent forms for stochastic programming problems can probably skip or skim this section.

Let us start by thinking about the two-stage stochastic problem. This is one in which an initial decision is made, time moves on, leaving the universe in some new state and then, in this state (i.e. after one time interval) a further decision, known as a *recourse decision*, is made. Such problems are known as *two-stage recourse problems*, since an initial decision is made followed later by a second, which allows corrective (recourse) action to be taken in light of the events that have been observed to occur. Figure 6 shows such a two-stage recourse problem.

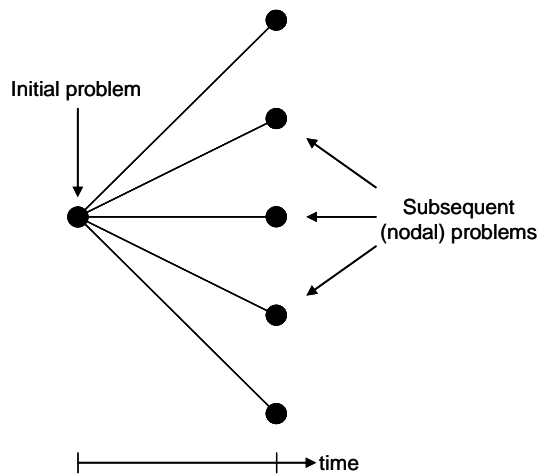

**Figure 6** *Example of a two-stage stochastic problem. Here there is uncertainty about which of the five possible successor states the world will move into after the first decision, so that decision must take that uncertainty into account.*

The solution to problems like this is taken to be the decision that maximizes the *expected* value of the objective function. This is not the only formulation we could use: for example we might instead want to maximize the objective value of the worst outcome, or of the best. Such formulations, however, lead to other types of problem, which we do not treat here (for example, robust optimization, which looks to ameliorate the worst case).

Given the aim of maximizing the expected objective function over all possible scenarios, the two-stage recourse problem can, in fact, be written as a single optimization problem. This formulation is known as the *deterministic equivalent*. An obvious solution method for (linear) stochastic programming problems (and one that is used in practice) is to convert the problem to its deterministic equivalent form and then solve it using the standard solution methods of linear programming (e.g. simplex or interior point methods, see e.g. Vanderbei (2000) for more details).

A common way of writing two stage problems with recourse is to do so as two linked optimization problems, with the second (recourse) problem dependent on both the solution chosen to the first and the evolution of the universe according to some particular scenario. We will follow this convention and then show that this can be rewritten easily as a single problem.

The first stage problem can be written in the following form

$$\min_{y} c'y + Q(y)$$
$$\text{s.t.} \quad Ay \leq b \quad (20)$$
$$y \geq 0.$$

Here $y$ is the first stage decision, $c'y$ gives the objective terms depending only on the first stage decision and $Q(y)$ represents the objective terms that come from the second stage recourse decisions. This $Q(y)$ term links the initial problem with the second stage recourse problems that we will see shortly and is dependent on the decision made in the first stage, since the recourse decisions will most likely depend on the initial decision.

In the second stage, we have an opportunity to make a further decision, knowing the outcome of whatever uncertain processes became certain in the time between the first and second decisions (for example, a consumption figure for an item during the problem period will be unknown at the first stage decision, but known for certain by the second). The problem we face in the second stage therefore depends on both the decision we made in the first ($y$) and the outcome of the uncertain processes, which we will call $u$. The second stage problem can thus be written as

$$[\,Q(y,u) = \,] \quad \min_{v} q(u)'v$$
$$\text{s.t.} \quad W(u)v \leq h(u) - T(u)y \quad (21)$$
$$v \geq 0.$$

This formulation looks a little gruesome, but really it is just a standard linear programming problem. However, all the terms that form the structure of the problem ($q$, $W$, $h$ and $T$) can depend on the information $u$ that has become available since the first decision. Splitting the right hand side of the constraint term into terms that depend on the initial decision and those that don't ($h$ and $T$, respectively) is done for convenience.

As previously mentioned, our problem formulation requires that we optimize the expected objective term given the possible evolution of the universe. This gives us a way of linking $Q(y)$ in (20) and $Q(y,u)$ in (21), since $Q(y)$ is simply the expectation of $Q(y,u)$ over all values of $u$ (i.e. over all possible evolutions of the universe from stage one to stage two).

$$Q(y) = \mathrm{E}_u\left[Q(y,u)\right] \qquad (22)$$

It is this relationship that would change if we wanted to use a different form of optimization (e.g. worst case, when the expectation would be replaced by min).

### A General Deterministic Equivalent Form

The following short section is rather unpleasant and can be skipped without difficulty.

The problems in (20) and (21) can be combined into a single problem, since both are minimizations. Since every possible realization of the future $u$ must be considered this formulation is a little ugly. The expectation over $u$ of (22) is here replaced by its integral form, and $(\Omega, P)$ is a *probability space* representing the future evolution of the universe (this is simply a way of denoting all the possible things that could happen $\Omega$ and assigning a probability to each of them and all reasonable combinations of them). The symbol $\omega$ represents one specific realization of the future in this universe and has some probability attached to it (though in the case of a continuous probability space this could be zero for any particular value of $\omega$). The sum (or integral in the case where the probability space is continuous) of their probability over all these must be equal to 1 (i.e. something must happen).

$$\begin{aligned}
\min_{y} \ & c'y + \int_{\Omega}\left(\min_{v_\omega} q(u_\omega)'v_\omega\right)dP \\
\text{s.t.} \quad & Ay \leq b \\
& W(u_\omega)v \leq h(u_\omega) - T(u_\omega)y && \forall \omega \in \Omega \\
& v_\omega \geq 0 && \forall \omega \in \Omega \\
& y \geq 0
\end{aligned} \qquad (23)$$

This problems looks, and is, difficult to solve. The difficulty comes from the unpleasant form of the expectation over all possible future realizations of the universe in the objective and from the two constraints that come from (21), which now have to be applied to every possible realization of the future. In the case of continuous probability spaces this leads to a continuous (and uncountably infinite) set of constraints. Furthermore, the optimization is over the continuous set of variables $y \cup \bigcup_{\omega \in \Omega} v_\omega$. This is not likely to be tractable in general.

**Discretization**

In order to make such optimization problems tractable in general, it is necessary to make an approximation to them. The key step in doing this is to make a discretization of the possible future evolution of the universe.

In many problems being considered, the universe will be able to evolve forward into any of an infinite number of possible states (for example, tomorrow's temperature in degrees Celsius could be any real number between, say, -50 and 80, of which there are an uncountably infinite amount). The distribution of this number might, however, be possible to estimate with some accuracy (tomorrow's temperature in Cambridge is predicted by the Met Office to be 8ºC and we might want to assume that the true temperature tomorrow will be normally distributed around this value with some variance, which we could perhaps discern from studying past weather forecasts and outcomes); unfortunately, it would be very surprising if it were 25ºC.

The idea of discretization is to approximate the probability distribution of the future with a finite number of *scenarios*, in each of which the universe moves into some specified state. When taken together, the (discrete) distribution of the unknown variable or variables across these scenarios should reflect the hypothesized (continuous) distribution of the variable(s) being modelled (in the case of multiple variables they should reflect the joint distribution of the variables). The nice thing about using such discretizations is that they make it possible to use all sorts of probability distributions, not just the obvious ones.

Returning to our temperature example, if we were using 3 scenarios, we might choose one in which the temperature tomorrow was 6ºC, one in which it was 8ºC and one in which it was 9ºC. As we add more scenarios, we want to the distribution of these numbers to more closely approximate the normal distribution that we hypothesise for tomorrow's temperature. Figure 7 shows the distribution of temperatures in the scenarios using 10, 100 and 1000 scenarios. The black line shows the normal distribution being approximated. As the number of scenarios improves, the approximation more closely approximates the hypothesized one. There are many sampling techniques to improve the distribution of these samples for small sample sizes, but these are a whole area of research on their own and are beyond the scope of this article. See, for example, Glasserman (2004) for examples of sampling methods for general Monte-Carlo techniques.

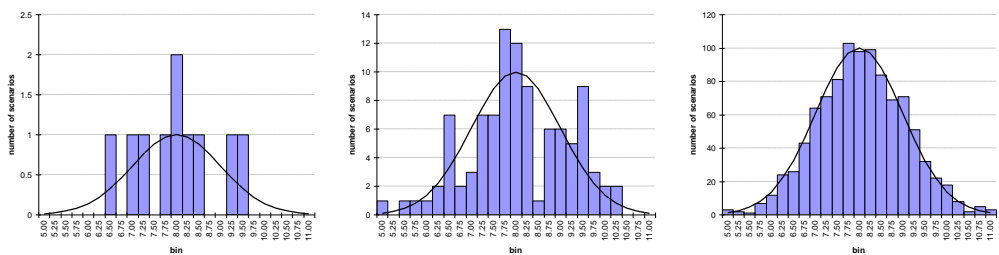

**Figure 7** *Distribution of temperature in individual scenarios with 10, 100 and 1000 scenarios, respectively*

Making such a discretization of the future state of the universe (represented as $u$) allows the expectation in (22) to be written in the form of a sum.

$$Q(y) = E_u[Q(y,u)] \triangleq \sum_{\omega \in \Omega} p(\omega) \cdot Q(y, u_\omega) \tag{24}$$

Here, $\omega$ represents a single discrete scenario, which occurs with probability $p(\omega)$. The state of the universe after one time period in scenario $\omega$ is given by $u_\omega$ and $\Omega$ is the (finite) set of all scenarios. The symbol $\triangleq$ is used here to denote "is discretized to".

Putting (24) into (20) and adding the (discretized) constraints from (21) we obtain the discretized deterministic equivalent form of the original problem.

$$\min_{y,v} \left\{ c'y + \sum_{\omega \in \Omega} p(\omega) \cdot q(u)' v_\omega \right\}$$

s.t.
$$\begin{aligned} Ay &\leq b \\ y &\geq 0 \\ W(u)v_\omega &\leq h(u) - T(u)y & \forall \omega \in \Omega \\ v_\omega &\geq 0 & \forall \omega \in \Omega \end{aligned} \tag{25}$$

The minimization can be performed over the variables from the first stage problem $y$ and the second stage problem $v$, where $v$ represents the set of all $v_\omega$ variables over all second stage scenarios $\omega$ (i.e. $v = \{v_\omega \mid \omega \in \Omega\}$).

The problem here is now a normal linear minimization problem, since the objective function is linear and the constraints are now finite in number as the number of scenarios is now finite. This problem is the *deterministic equivalent form* of the discrete stochastic programming problem, where the future is represented by one of a finite number of successor states to the present. Although this discretization might only be an approximation of the real situation, it is important because it can be solved using the standard methods of linear programming.

This is, in fact, the problem that we are interested in solving in this article by applying Benders decomposition. In light of our earlier discussion of Benders, this problem looks like a promising candidate for decomposition, since the problem is easily divided into the first stage problem and a number of second stage problems, with one coming from each of the future scenarios being considered. The constraints in (25) can all be ascribed to exactly one of these problems. Furthermore, the second-stage problems are independent of each other (that is, they do not share variable or constraints) and are only linked by the choice of the first stage variables, $y$. This is ideal, since the only information from the rest of the problem that needs to be passed to each of the sub-problems is the first stage decision. Once this is made the second stage problems can be solved independently. This property means that a solution method based on Benders decomposition lends itself to parallelization.

In our first section on Benders we examined a decomposition into only two sub-problems, but this can be readily extended to multiple sub-problems, as we shall see in the next section.

## Applying Benders to Stochastic Programming Problems

The preceding section illustrated how stochastic programming problems have a structure that lends itself naturally to decomposition. In this section we will see how that decomposition leads to an elegant algorithm for the solution of stochastic programming problems. The application of decomposition methods to stochastic programming problems (as opposed to mixed integer problems, which is what Benders dealt with) was first done by Van Slyke and Wets (1969) and is sometimes known as the *L-shaped method*. Rico-Ramirez (2002) gives a readable tutorial on this material.

The first step in applying Benders decomposition is knowing what problem we are trying to solve. In this case it is the deterministic equivalent form of the stochastic programming problem that we are interested in, given by (25). This becomes our problem P from the *Basic Ideas* section. The decomposition we use is into a single first stage problem and a number of second stage problems, one corresponding to each of the possible scenarios that could be realized in stage 2.

From the *Basic Ideas* section we saw that the problem $P^k$ was equivalent to P, but used constraints derived from the solutions of the dual of the sub-problems. If, as in the case of stochastic programming problems, we have many sub-problems then $P^k$ will contain constraints deriving from the duals of each of these sub-problems. We also saw that the idea of Benders decomposition as a solution method was to solve the simpler problem $P^r$ (a form of $P^k$ only incorporating $r$ of $P^k$'s $k$ constraints) and hope that this would give us the optimal solution before we had added back all $k$ constraints. The problem that was not resolved in that section was how to add the constraints to the problem $P^r$ as $r$ was increased in a way that would make this more likely. In this section we will address this issue.

The first step in any solution method is to come up with an initial estimate of the first-stage variables $y$. This can be arbitrary. We then use this value of $y$ as a fixed quantity and attempt to solve each of the sub-problems assuming that this chosen value was the decision made in the first stage. There are then two possibilities: either some of the sub-problems are infeasible with this chosen value of $y$, or all the sub-problems are feasible. Both of these situations give us a way to add constraints to the problem. In the first situation we will be provided with constraints known as *feasibility cuts* and in the second we will obtain constraints known as *optimality cuts*. We will examine each in turn.

### Feasibility Cuts

If some of the sub-problems are infeasible with the chosen initial decision $y$, it tells us that that value of $y$ is not feasible for the entire problem (since for a solution to be feasible the initial decision must lead to a feasible second stage problem no matter what the future evolution of the universe). We want to use this information to add a constraint back into the reduced problem that we are solving P$^r$. Such constraints are known as feasibility cuts.

An example might make the situation clearer. Consider the case where the initial decision $y$ is a two dimensional vector and where each of the sub-problems impose a single feasibility constraint on the value of $y$. That is, one of the sub-problems will become infeasible when $y$ violates the corresponding constraint from that sub-problem. Graphically, we could illustrate the problem as follows.

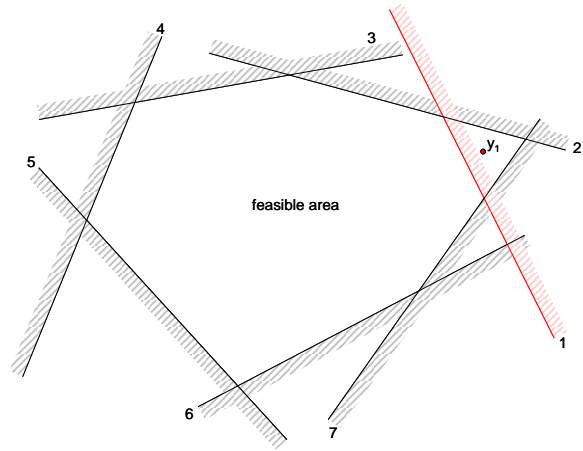

**Figure 8** *Step 1: Feasibility constraints for the example problem. The initial guess is shown as $y_1$. Violated constraints are shown in red.*

Here, the initial choice of $y$ is given by the point $y_1$. The problem shown has seven constraints that must be satisfied to make all the sub-problems feasible. These could come from any number of sub-problems (including more than seven, since some sub-problems might be feasible for all choices of $y$), but here we will consider the case where one constraint comes from each of seven sub-problems. From the point of view of a stochastic programme, this implies a two-stage structure with an initial decision and a universe that can evolve into one of seven subsequent states at the second stage. In the situation shown in Figure 8 the initial guess $y_1$ causes sub-problem 1 to become infeasible. Therefore, the constraint arising from this sub-problem is added to the set of problem constraints on $y$ and the problem is re-solved for $y$ (the way in which such constraints are derived has not yet been described, but is covered below). The new situation is shown in Figure 9.

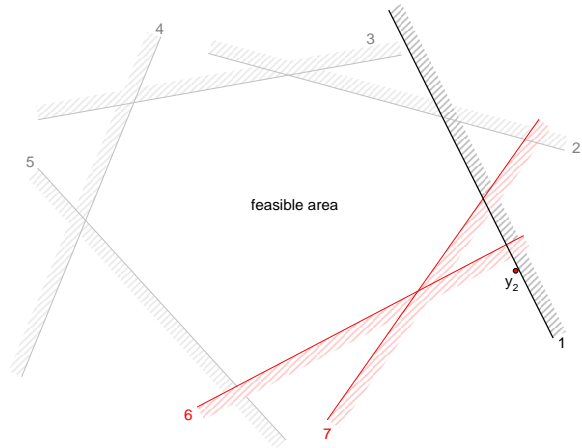

**Figure 9** *Step 2: Constraints included in the problem solved at this step are shown in black. Constraints violated by the new proposed solution $y_2$ are shown in red.*

The new solution $y_2$ satisfies the constraint that ensures the feasibility of sub-problem 1, but now violates those that ensure feasibility of sub-problems 6 and 7. These constraints are now added to the problem in order to find a new proposed solution. This is shown in Figure 10. (In fact, if the sub-problems are solved sequentially it can also be the case that a single feasibility constraint is added as soon as any sub-problem is found to be infeasible. Both methods will work; the first is a sort of multicut, where several cuts can be generated at each step of solution and the latter is a sort of unicut, where only a single cut is generated at each step of the solution. We will return to unicut and multicut in the section on *Optimality Cuts*, below).

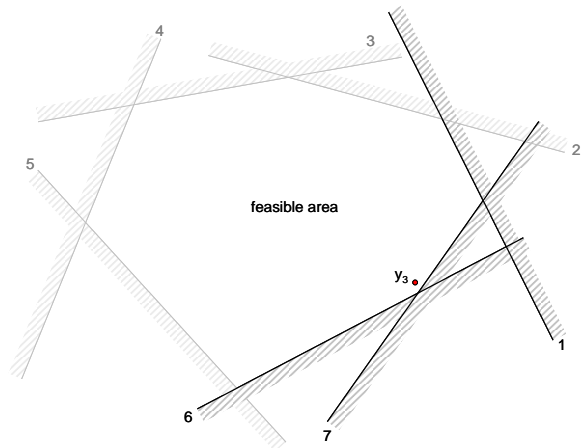

**Figure 10** *Step 3: Constraints ensuring the feasibility of sub-problems 6 and 7 are added to the problem and a new proposed $y_3$ is found.*

The new proposed solution $y_3$ is found by solving the optimization for *y* including constraints 1, 6 and 7. This proposed solution does not violate any of the feasibility constraints from the sub-problems and so this choice of the first stage decision *y* results in all second-stage problems being feasible. Furthermore by only adding constraints deriving from sub-problems that have been made infeasible by the proposed solution at each step, we have found this point using only three of the seven sub-problem feasibility constraints. However, this might still not be the optimal choice of *y* (merely a feasible one) and we might need to add a series of further constraints known as optimality cuts to help us find this optimal point (covered in the next section). The remainder of this section will be devoted to deriving the feasibility cuts that we have been discussing.

In the above example, we assumed that we could somehow derive the constraints on the initial problem from the sub-problems that ensure feasibility. Happily, this is the case and can be done by considering a modified version of the sub-problem and its dual.

The form of the sub-problems is given by (21), which we will rewrite here without the dependencies on *u* for notational convenience.

$$[\,Q(y,u) = \,] \quad \min_{v} q'v$$
$$\text{s.t.} \quad Wv \leq h - Ty \quad (26)$$
$$v \geq 0$$

For feasibility the value of the objective function is not relevant; all that matters is that the constraints are satisfied. We can adapt the problem in (26) to test whether this is the case in such a way that it will always be feasible. This is done by adding in a pair (of vectors) of unconstrained (or *slack*) variables, $v^+$ and $v^-$. We need to introduce one free variable (i.e. one $v^+$, $v^-$ pair) for each constraint in problem (26), so that each constraint can always be satisfied.

$$[\,V = \,] \quad \min_{v,v^+,v^-} e'(v^+ + v^-)$$
$$\text{s.t.} \quad Wv + v^+ - v^- \leq h - Ty$$
$$v \geq 0 \quad (27)$$
$$v^+ \geq 0$$
$$v^- \geq 0$$

Here *e* is a vector of 1s of the appropriate dimension (i.e. the number of constraints in (26)). The variables $v^+$ and $v^-$ ensure that the constraints from (26) can always be satisfied, since, as both are unbounded below, the term $v^+ - v^-$ is unbounded in either direction and can thus always be chosen to satisfy the constraints. The objective here measures (to be precise, measures in the 1-norm) the amount by which the constraints are violated, since it is just a sum of the values assigned to $v^+$ and $v^-$. We therefore want to minimize this violation of the constraints and if we can reduce it to zero it means that the constraints in the original sub-problem (26) can be satisfied, implying that it is feasible.

If, however, on solving (27) we find that we cannot reduce the objective to zero it means that the sub-problem is infeasible. In this case we can add a feasibility cut, designed to add a constraint back to the master problem that will ensure that this sub-problem will be feasible. Remember, though, this process only reveals feasibility constraints that were violated with the current value of *y*, so when we re-solve the master problem with these new feasibility cuts we will get a new optimal value of *y* and we will have to repeat this process.

To devise a constraint that we can use in the master problem we must come up with a constraint in terms of $y$ and some known quantities. We can do this by considering the dual problem of (27):

$$[V =] \quad \max_{\sigma} \sigma'(h - Ty)$$
$$\text{s.t.} \quad \sigma'W \leq 0 \qquad (28)$$
$$|\sigma| \leq e.$$

For the original sub-problem (26) to be feasible the optimal value of the modified sub problem (27) $V$ has to be zero. By the strong duality theorem (see *Aside on Duality Theory*) this means that if the optimal value of the dual of this modified problem (28) is zero then this will be the case and the sub-problem (26) will be feasible. If we solve this dual problem for a particular value of $y$, call it $y_v$ we will obtain a particular solution to the dual $\sigma_v$. If $y_v$ is such that it makes (26) infeasible this solution will result in a value of $V > 0$. We want to avoid this and so can set a constraint on $y$ (in the original first stage problem) that makes this impossible. This is very simple to arrange, and can be done by placing a constraint on $V$ that requires it to be less than or equal to zero for the particular solution of the dual problem $\sigma_v$. This can be ensured by constraining the value of the objective term from (28) to be less than or equal to zero. The maximization in (28) will then always ensure that it is zero and that (26) is feasible. The constraint we need to add into the original problem is therefore

$$\sigma_v'(h - Ty) \leq 0. \qquad (29)$$

This is a constraint we can actually add to the master problem because it is a constraint in $y$ (the master problem solution) only and all the quantities $\sigma_v$, $h$ and $T$ are known.

Writing $D_v = \sigma_v'T_k$ and $d_j = \sigma_v'h$ we can rewrite (29) as

$$D_v y \geq d_v \qquad (30)$$

Note that this constraint will only ensure feasibility when the solution of the dual is $\sigma_v$. If a different value of $y$ was chosen it is possible that the dual (28) will have a different solution $\sigma_w$. If this new value of $y$ makes the sub-problem (26) infeasible again then it will be necessary to add another feasibility constraint to the first phase problem, using the dual solution $\sigma_w$ in place of $\sigma_v$ in (29).

However, since (27) is always feasible, its dual is always bounded and can only have one of a finite number of solutions (corresponding to the corner points of its feasible region). This means that the addition of such feasibility constraints must be a finite process.

What does it mean to say that the constraint will only ensure feasibility when the dual has some particular solution $\sigma_v$? Well, the dual problem is an optimization problem with as many variables as there are constraints in the primal problem. It can therefore have a number of cornerpoints to its feasible region and thus a number of candidate solutions. Since the constraints of the primal problem and thus the objective of the dual can depend on the value of the solution from the first-stage parent problem $y$, the optimal value of the variables of the dual solution can also depend on $y$. The feasibility constraint generated is the one corresponding to the particular solution $y$ of the first-stage problem (the particular solution of the dual for the given value of $y$ for the $v$th version of the parent problem is $\sigma_v$). By the nature of optimal solutions being at cornerpoints of the feasible region (see *Aside on Linear Programming*) it is often the case that the dual solution remains the same for some range of first-stage solutions $y$, so what we are really saying is that the feasibility constraint imposed will be relevant to the solution of the first-stage problem (imposing feasibility for the given sub-problem) for all values in this range of values of $y$. The example in the next section might make this clearer.

The process we have described in this section gives us a way of generating constraints for the original problem when the sub-problems prove infeasible. However, it is possible that none of the sub-problems will be infeasible without having found an optimal solution. To see this, consider a stochastic programme with full recourse where whatever decision is made in the first stage, the second stage problems always have a solution. In such a problem there will not be any feasibility cuts and so we need another way of generating constraints using information from the sub-problems in the case when they are all feasible. Such is the case in, for example, problems with *full recourse* where there is always some subsequent course of action that will achieve a feasible solution, no matter how stupid the initial decision. This leads us to the next topic, *optimality cuts*, to which we will turn after presenting a simple example of generating a feasibility cut.

### Example of Feasibility Cut

This section will present an example of a feasibility cut for a simple problem, which will hopefully illustrate some of the ideas of the previous section in a more concrete way. The example in this section is a maximization rather than a minimization as shown above but that does not affect the point being illustrated (and one type of problem can be transformed into the other by negating the objective term).

Consider a two stage optimization problem with some first stage problem with an initial decision $y$ (a single variable) and a sub-problem with a decision $v$ (again a single variable). The exact form of the first stage problem does not matter in this case: we just assume that it gives us a decision $y$. The sub-problem SP is

$$\text{SP:} \quad \begin{aligned} \max \quad & v \\ \text{s.t.} \quad & v \geq 1 \\ & v \leq 3 - y \\ & v \leq y + 3. \end{aligned} \qquad (31)$$

Note that here the constraints in the sub-problem depend on the initial decision $y$. When considering the sub-problem, we can simply regard this initial decision as a parameter of the problem. Figure 11 shows the problem. Clearly, the problem is only feasible for values of $y$ between -2 and 2. Outside this range there is no value of $v$ that can be chosen so as not to violate at least one of the constraints.

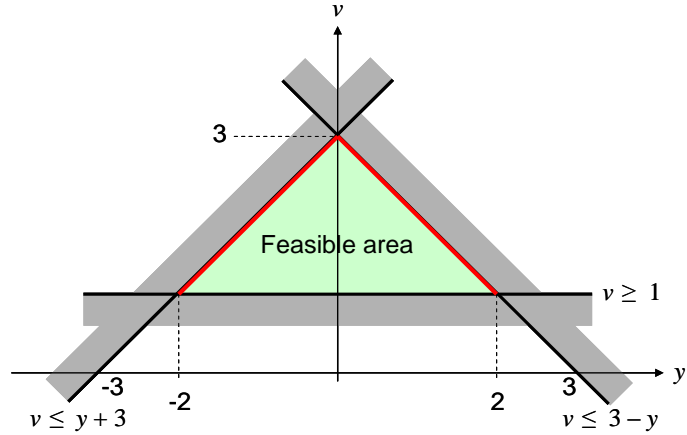

**Figure 11** *Illustration of problem (31). The red lines show the optimal values of v for each particular value of y. The problem is only feasible for $-2 \leq y \leq 2$*

The first step in the feasibility cut part of the solution process is forming the feasibility problem for the sub-problem (corresponding to (27) above). For this problem the feasibility problem FP is

$$
\text{FP:} \quad
\begin{aligned}
\min \quad & v_1^+ + v_1^- + v_2^+ + v_2^- + v_3^+ + v_3^- \\
\text{s.t.} \quad & v + v_1^+ - v_1^- \geq 1 && (\alpha) \\
& v + v_2^+ - v_2^- \leq 3 - y \quad \equiv \quad -v - v_2^+ + v_2^- \geq y - 3 && (\beta) \\
& v + v_3^+ - v_3^- \leq y + 3 \quad \equiv \quad -v - v_3^+ + v_3^- \geq -y - 3 && (\gamma) \\
& v_i^+ \geq 0 \quad && i = 1, 2, 3 \\
& v_i^- \geq 0 \quad && i = 1, 2, 3.
\end{aligned}
\quad (32)
$$

All we have done here is to add one pair of slack variables $v_i^+$ and $v_i^-$ to each constraint in the sub-problem and replaced the objective by one that minimizes the feasibility violation of the problem. For the constraints labelled ($\beta$) and ($\gamma$) the equivalent greater than constraint is shown as this puts the minimization problem into a standard form and makes formation of the dual a little more intuitive.

If we solve the feasibility problem and find that its optimal objective value is greater than 0 (i.e. there is some feasibility violation) then we want to form a feasibility cut. The next step is then to form the dual of the feasibility problem. Since (32) is a minimization what we are looking for in forming the dual is a lower bound on the objective value. Dual formation is therefore done by taking multiples of the constraints in such a way that we come up with an expression to which the objective is greater or equal. This is essentially the same process as that described in the Aside on Duality Theory, but with the direction of the inequalities reversed. The dual of the feasibility problem is

DFP:
$$\begin{aligned}
\max \quad & \alpha + \beta(y-3) + \gamma(3-y) \\
\text{s.t.} \quad & \alpha - \beta - \gamma \leq 0 & (v) \\
& \alpha \leq 1 & (v_1^+) \\
& -\alpha \leq 1 & (v_1^-) \\
& -\beta \leq 1 & (v_2^+) \\
& \beta \leq 1 & (v_2^-) \\
& -\gamma \leq 1 & (v_3^+) \\
& \gamma \leq 1 & (v_3^-) \\
& \alpha, \beta, \gamma \geq 0.
\end{aligned}$$
(33)

The variables from the feasibility problem to which each of the constraints apply are shown in brackets. This dual is a linear optimization problem in three variables. We can visualize its feasible area as a volume in the space spanned by these three variables and, remembering that solutions to linear optimization problems lie on the cornerpoints of the feasible area, we can enumerate the list of possible solutions to the problem. Figure 12 shows the feasible area of problem (33).

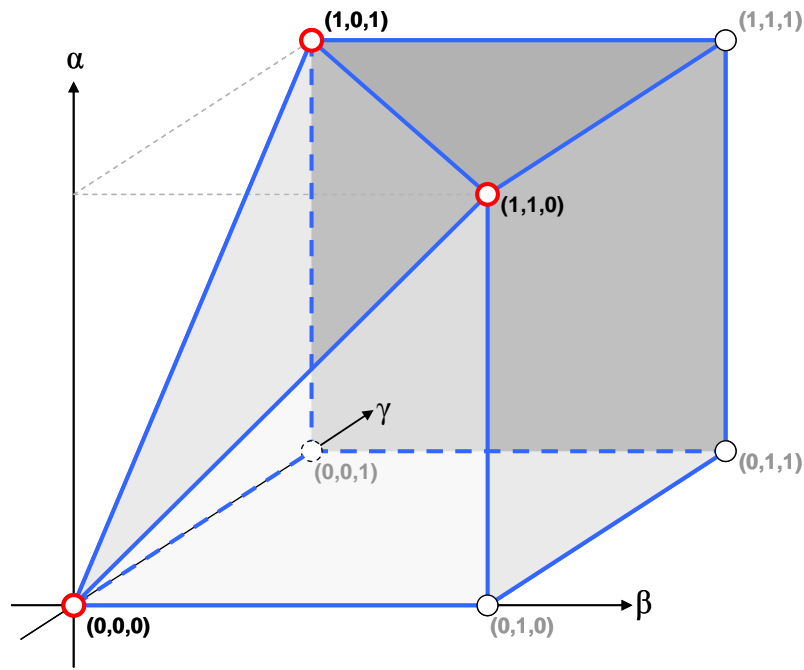

**Figure 12** *Feasible region for problem (33), outlined in blue. Cornerpoints are shown as circles, with optimal solutions shown as red circles.*

The list of potential solutions is given in the following table.

| Solution ($\alpha, \beta, \gamma$) | Objective value | Note |
|---|---|---|
| (0, 0, 0) | 0 | A |
| (0, 0, 1) | $-y-3$ | dominated by B |
| (0, 1, 0) | $y-3$ | dominated by C |
| (0, 1, 1) | -6 | dominated by A |
| (1, 0, 1) | $-y-2$ | B |
| (1, 1, 0) | $y-2$ | C |
| (1, 1, 1) | -5 | dominated by A |

We can see that we only really need to consider three of these potential solutions, A, B and C, as the others are all dominated by them (i.e. the solution at one of A, B or C will always be greater, no matter what the value of $y$). In fact, it is fairly obvious that A will be the solution when $y$ is between -2 and 2, B will be the solution when $y$ is less than -2 and C will be the solution when $y$ is greater than 2.

We can now use these solutions to the dual problem to devise feasibility cuts. Now, the primal problem is feasible when $y$ is between -2 and 2, so in this case we will be trying to generate a feasibility cut. Solution A, therefore, does not give us a feasibility cut. However if $y$ is less than -2 or greater than 2 then the sub-problem was infeasible and we will try to form a feasibility cut. Remember that we can do this by insisting that the objective of the dual is smaller or equal to zero. This gives the following feasibility cuts:

| Condition | Dual solution ($\alpha, \beta, \gamma$) | Feasibility cut |
|---|---|---|
| $y < -2$ | (1, 0, 1) | $1+1\times(-y-3) \leq 0 \equiv $ **$y \geq -2$** |
| $y > 2$ | (1, 1, 0) | $1+1\times(y-3) \leq 0 \equiv $ **$y \leq 2$** |

So here are our eminently sensible feasibility cuts. If $y$ is less than -2, we get a feasibility cut that insists $y$ be greater than or equal to -2; if $y$ is greater than 2 we get a feasibility cut that insists that $y$ is less than or equal to 2.

It should now be much clearer what we were talking about at the end of the previous section when we said that the feasibility cut ensures feasibility for a particular solution of the dual problem. Each of feasibility cut corresponds to a particular solution of the dual of the feasibility problem and any value of the first stage decision y that leads to the same solution of the dual (here, for example, all values of y greater than 2) will generate the same feasibility cut, which will ensure that infeasible values of the first stage decision corresponding to that particular solution of the dual cannot be chosen in the future.

## Optimality Cuts

After generating feasibility cuts an optimal solution probably still will not have been obtained. In order to reach an optimal solution rather than a merely feasible one we also need to generate optimality cuts, which send information back to the first stage problem about how to make the first-stage decision in such a way as to be optimal. Generating these cuts is the focus of this section.

When considering generating optimality cuts we may start with the assumption that all sub-problems are feasible with the chosen value of the first-stage problem $y$, since if they were not then feasibility cuts would have been generated.

Let us re-examine the original optimization problem:

$$\min_y \; c'y + Q(y)$$
$$\text{s.t.} \quad Ay \leq b \tag{34}$$
$$y \geq 0.$$

As described above, $Q(y)$ is the expected value of the bits of the objective function that come from the sub-problems given that $y$ is the choice of the first-stage decision variables.

$$Q(y) = \mathrm{E}_\omega Q(y, u_\omega) \tag{35}$$

Here $\omega$ is one of the successor states of the root node. Each one of these successor states has a sub-problem associated with it (as mentioned above in *Stochastic Programming Problem Structure*) and in each state $u_\omega$ is the realization of the random processes corresponding to that state. If there are only a finite number of discrete successor states then we can write (35) as

$$Q(y) = \sum_\omega p_\omega Q(y, u_\omega), \tag{36}$$

where $p_\omega$ is the probability of state $\omega$ occurring. In this case $p_\omega Q(y, u_\omega)$ is the part of the overall objective function corresponding to the sub-problem associated with state $\omega$.

The function $Q(y)$ is a complex and expensive function to evaluate because it depends on finding an optimal solution to all of the sub-problems. Furthermore, considered simply as a problem in $y$, the optimization problem (34) has a non-linear objective function, since $Q(y)$ is non-linear in $y$ (we will return to the nature of $Q(y)$ in a moment).

Remember that what we are looking for is a sensible way to add constraints to some master problem, a less constrained version of the complete problem, so that we can come up with sensible estimates of $y$ without solving the entire problem. In Benders decomposition we start with the first stage problem as the initial master problem and the constraints that we add to it correspond to feasibility and optimality of the (second-stage) sub-problems. We hope that by adding constraints sensibly our estimates will converge to the optimal solution whilst only requiring us to solve a much smaller problem.

So here, we try the following:

- Remove the complicated $Q(y)$ term from the objective function of (34)
- Replace it with a variable $\theta$
- Use information from the sub-problems to place bounds on $\theta$, which we call optimality cuts.

The idea is that, though eventually these bounds will end up recreating the function $Q(y)$, at first they will not, and the problem will be much simpler. Hopefully we will find the optimal solution long before the original $Q(y)$ function is replicated.

The first thing we need to understand is what this $Q(y)$ function is actually like.

Since $Q(y)$ is a weighted sum of the objective functions of the sub-problems $Q(y,u)$ over each possible value of $u$, understanding the nature of these is a good first step. These functions are the optimal objective values of the second-stage sub-problems (where $v$ is the second-stage decision):

$$
\begin{aligned}
Q(y,u) = \quad & \min_{v} q(u)'v \\
\text{s.t.} \quad & W(u)v \leq h(u) - T(u)y \\
& v \geq 0
\end{aligned}
\tag{37}
$$

Perhaps the easiest way to understand what these are like is with a very simple example. Consider the following second-stage problem (where, for simplicity, there is no dependence on $u$).

$$
\begin{aligned}
Q(y,u) = \quad & \min_{v} v \\
\text{s.t.} \quad & v \geq y \\
& v \geq 10 - y \\
& v \geq 9 - 0.5y
\end{aligned}
\tag{38}
$$

This corresponds to the problem shown in Figure 13.

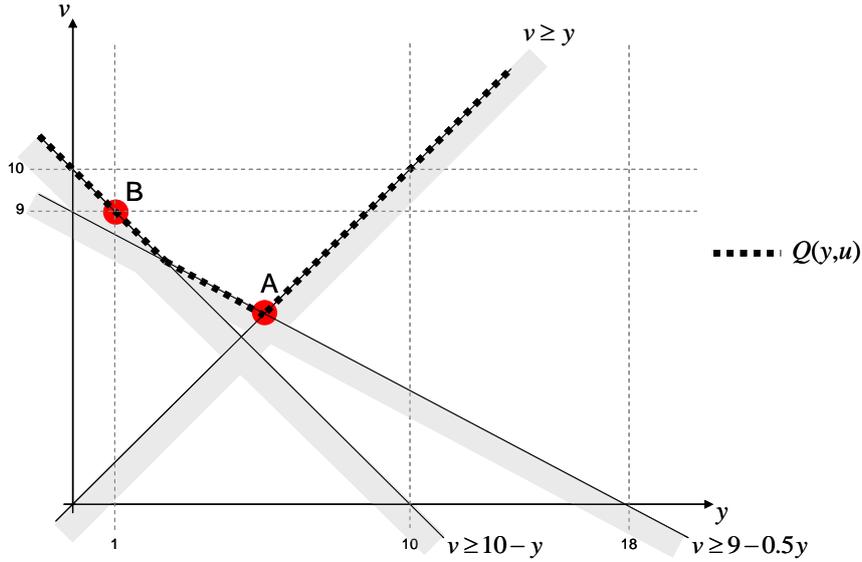

**Figure 13** *Example sub-problem corresponding to* (38). *Constraints are shown in grey. The objective function Q(y,u) as a function of y is shown as a heavy dotted line.*

The value that minimizes the objective function of this sub-problem is located at point A. However, in the sub-problems we do not get to choose the value of *y*; that comes to us from the first stage. If, say, the value of *y* from the first stage was 1 then the optimal value of the objective is at point B. For each sub-problem the value of *u* is fixed to $u_\omega$, so even if we had dependency on *u* (which we don't in this simple example) for any given sub-problem $Q(y, u_\omega)$ will really just be a function of *y*. Figure 13 illustrates the type of function that it will be: piecewise linear and convex. It will always be piecewise linear because it is constructed from a finite number of linear constraints. To see that it will always be convex, consider adding any further constraint to the function. It can only make it 'more' convex since any additional constraint will cause the value of the objective to be greater or equal to the value of the objective without the constraint at every point.

Now, *Q(y)* is the expectation over all second-stage states *ω* of these sub-problem *Q(y,u)*. When we have a finite number of successor states this is a finite sum, as in (36). Since the sum of two convex functions is itself convex then *Q(y)* must also be a convex function in *y*. It is also piecewise linear because it is a sum of other piecewise linear functions. This is sufficient knowledge of *Q(y)* to generate constraints for the master problem, as we now will demonstrate.

To generate further master problem constraints (known as *optimality cuts*) we replace the complicated *Q(y)* term in the objective of (34) by a single variable *θ* and then attempt to bound this. To bound it we use the sub-problems.

Consider the sub-problem for the case when the random processes *u* take the realization $u_\omega$. Call this sub-problem *ω*.

$$Q(y, u_\omega) = \min_v q(u_\omega)'v$$
$$\text{s.t.} \quad W(u_\omega)v \leq h(u_\omega) - T(u_\omega)y \quad (39)$$
$$v \geq 0$$

If we take the dual of this we get

$$Q(y, u_\omega) = \max_\pi \pi'\big(h(u_\omega) - T(u_\omega)y\big)$$
$$\text{s.t.} \quad \pi'w(u_\omega) \leq q(u_\omega)' \quad (40)$$
$$\pi \geq 0$$

Let's say that $y$ (the first stage decision) is that obtained from the $i^{th}$ iteration of the master problem and call it $y^i$. Now let $\pi^i_\omega$ be the optimal values of the $\pi$ variables in (40), the dual of the sub-problem (39) corresponding to $u_\omega$ and $y^i$, so that (dropping the dependence on $u_\omega$ for notational convenience)

$$Q(y^i, u_\omega) = (\pi^i_\omega)'\big(h - Ty^i\big). \quad (41)$$

This expression comes from the objective of (40).

We have, therefore, an expression for $Q(y^i, u_\omega)$ at $y^i$ in terms of $y^i$. Now, because we know something about the shape of $Q(y, u_\omega)$ we can use this to infer something more general. $Q(y, u_\omega)$ is piecewise linear and convex, so, in the region of $y^i$ on the same linear piece the value of $Q(y, u_\omega)$ just varies linearly with $y$ and we have

$$Q(y^i + \Delta, u_\omega) = (\pi^i_\omega)'\big(h - T(y^i + \Delta)\big). \quad (42)$$

Figure 14 shows the situation.

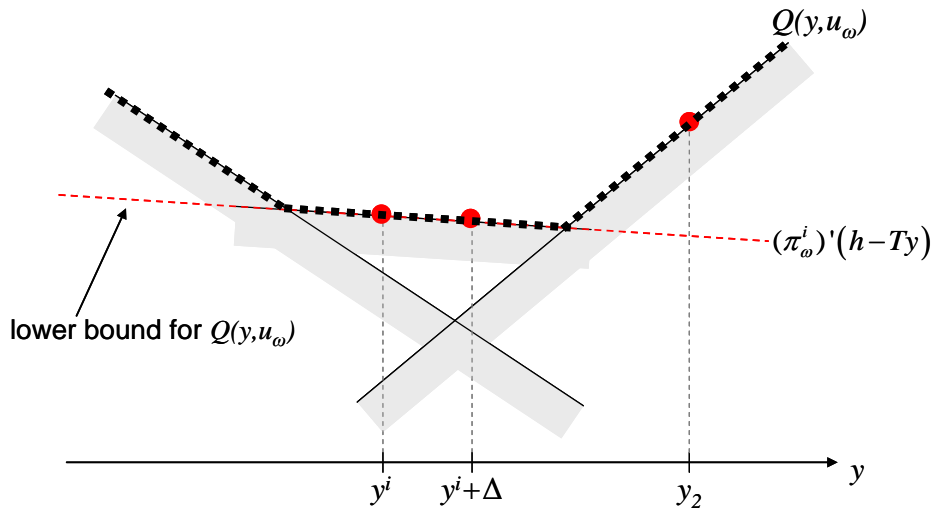

**Figure 14** *Example of creating a lower bound for a sub-problem objective function Q(y,u) using the dual objective term evaluated for some particular initial decision $y^i$.*

If we move further away, onto a different linear piece (e.g. $y_2$ in Figure 14), we know by the convexity of $Q(y, u_\omega)$ that

$$Q(y_2, u_\omega) \geq (\pi_\omega^i)'(h - Ty_2) \tag{43}$$

How do we know that the right hand side of (41) corresponds to the linear piece that forms part of $Q(y, u_\omega)$ at $y^i$? The answer comes from duality theory. The linear piece of $Q(y, u_\omega)$ is produced by a constraint in the primal sub-problem. Recall that when we form the dual each constraint in the primal problem leads to a variable and each primal variable leads to a dual constraint, and the constraints from the primal problem go into the objective of the dual problem multiplied by one of the dual variables (which can be thought of as Lagrange multipliers if you are so inclined). The $\pi'(h - Ty)$ in the objective of the dual of sub-problem (40) correspond directly to the constraints in the primal. Since it is these primal constraints that form $Q(y, u_\omega)$, this means that the $\pi'(h - Ty)$ terms really do correspond to the linear pieces of $Q(y, u_\omega)$.

What we are effectively doing with (43) is using constraints from the sub-problem to produce a lower bound on its objective $Q(y, u_\omega)$. By using the dual problem we are able to express this in terms of the initial decision $y$. The constraints that we are using to produce this bound, as shown in Figure 14, are the constraints that are biting at the point $y^i$. (The constraints that are in effect are determined by the $\pi_\omega^i$ vector).

Now, since $Q(y)$ is simply a weighted sum of these sub-problem objectives as shown in (36), we can use the bound we have on the individual sub-problem objectives to create a bound on $Q(y)$ by putting the bounds from (43) into the expression for $Q(y)$.

$$Q(y) \geq \sum_\omega p_\omega (\pi_\omega^i)'(h - Ty) \tag{44}$$

Remember that we were trying to create a lower bound on $Q(y)$ which we were then going to apply to $\theta$. Using (44) we can now do this by insisting that

$$\theta \geq \sum_\omega p_\omega (\pi_\omega^i)'(h - Ty). \tag{45}$$

Letting $g = \sum_\omega p_\omega (\pi_\omega^i)' h$, and $G = \sum_\omega p_\omega (\pi_\omega^i)' T$, we can write this as

$$\theta \geq g - Gy \tag{46}$$

This is a bound on $\theta$ in terms of $y$ that we can add to the master problem. We have arrived at our optimality cut for the $(i + 1)^{th}$ iteration.

**Termination**

At this point we can also test if we have reached an optimal solution. If the value we have previously calculated for $\theta$ in the master problem already satisfies our new optimality cut (46) then we have, in fact, found an optimal solution. Why is this? Well, if the initial decision was chosen to be optimal in light of the current approximation to the sub-problem objective function (made up of existing optimality cuts) then if, at the point chosen as optimal for the first stage, this approximation turns out to be accurate, the initial decision is optimal in light of the true sub-problem objective and so is the overall optimal. We can test this by seeing whether the value of our approximation $\theta$ is already equal to (or greater than) the value of the currently active constraint on it, i.e. the optimality cut in (46). If this is the case then we have reached the true solution and we can terminate. The test for termination is therefore

$$\theta^i \geq g - Gy, \tag{47}$$

where $\theta^i$ is the value of $\theta$ from the most recent solution of the master problem (i.e. the value of $\theta$ on the $i$th iteration of the algorithm). If (47) is satisfied, then we terminate. Figures 16 and 17 in the illustration below give a graphical example of this and why it works.

If, on the other hand, we do generate a further optimality cut we must iterate the entire algorithm again, including the feasibility tests that generate feasibility cuts, as adding in a new optimality cut can lead to new infeasibilities.

**Illustration of Optimality Cuts**

It is sometimes said that optimality cuts create an outer linearization of $Q(y)$. This simply refers to the fact that the optimality cuts, as we saw above, place linear lower bounds on the value of $Q(y)$. As more optimality cuts are added these will approximate the shape of $Q(y)$ more and more closely, but always remain below (or 'outside') it at some points until all the constraints have been added.

We will attempt to illustrate this with an example. For simplicity consider a problem with only one sub problem (so, in this case $Q(y)$ will coincide with $Q(y, u_\omega)$ for the single sub-problem). Say the first solution to the master problem is $y_1$. This is passed to the sub-problem and we then minimize the sub-problem based on this value of $y$.

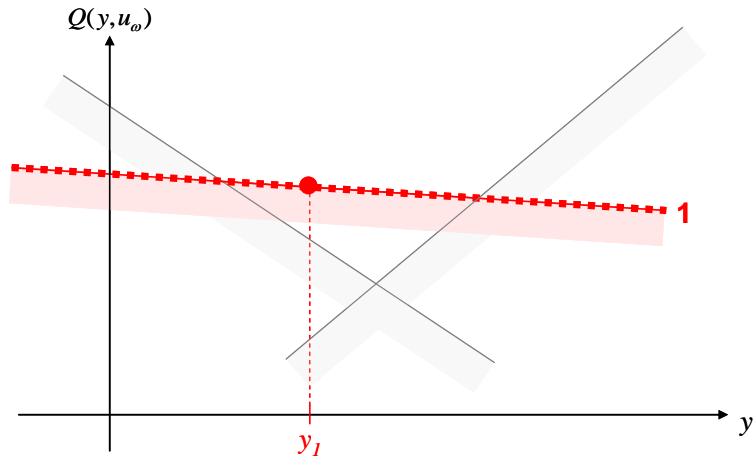

**Figure 15** *Step 1: Minimize the sub-problem with $y_1$ as the value of the initial decision. The active constraint is shown in red. The heavy dotted line shows the current approximation to $Q(y, u_\omega)$.*

Constraint 1 (shown in red in Figure 15) is the active constraint at point $y_1$. So, as an optimality cut we pass back the fact that we know $Q(y)$ (the sub-problem objective in this special case) lies above this. After adding this new constraint to the master problem and solving again we obtain a different initial decision, $y_2$. This is again passed to the sub-problem, which is re-solved.

This initial decision $y_2$ is the optimal initial decision if the actual objective of the sub-problem is only constrained by the optimality constraints added thus far (i.e. constraint 1) at the point $y_2$. The blue circle in figure 8 shows the objective value that the optimizer 'thinks' it will get from the sub-problem given its current constraints on that value. But, instead of getting point A as the objective of the sub-problems the optimizer actually gets point B since at point $y_2$ constraint 2 (which it did not know about when choosing point $y_2$) is active.

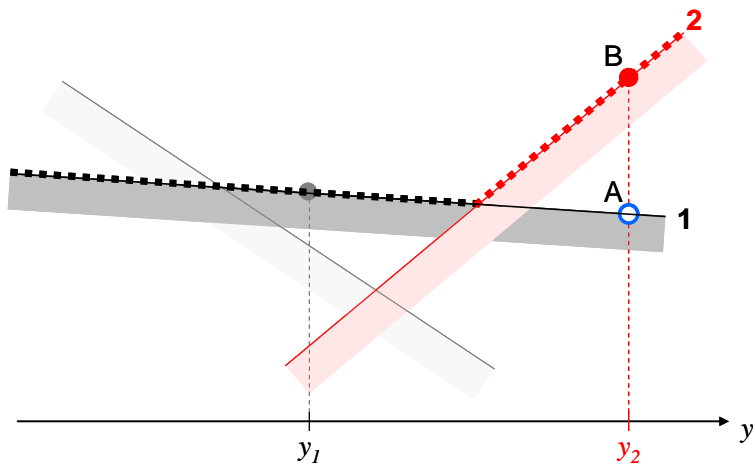

**Figure 16** *Step 2: Minimize the sub-problem with $y_2$ as the value of the initial decision. Another piece is added to the approximation of $Q(y, u_\omega)$.*

Constraint 2 is now the active constraint and we know that $Q(y)$ must also lie above this constraint. This is then passed back to the master problem as a further optimality cut. The heavy dotted lines show the lower bounds applied to $\theta$ at each step. As we add more constraints through optimality cuts this approximation of $Q(y)$ by $\theta$ will look increasingly similar to $Q(y)$. It will always lie below it at some points (until all the constraints have been added) and so is an 'outer' approximation of $Q(y)$.

Again we solve the master problem with the new optimality cut added to it. We get the solution $y_3$. (In this case we would expect this solution to lie to the left of $y_2$ because even when the optimizer thought $Q$ would be lower going right, when it did not have the second constraint in place for $Q$, it chose $y_2$ as the optimal solution. This combined with the fact that the part of the objective function from the primal solution is convex in $y$ will mean that $y_3$ should lie to the left of $y_2$).

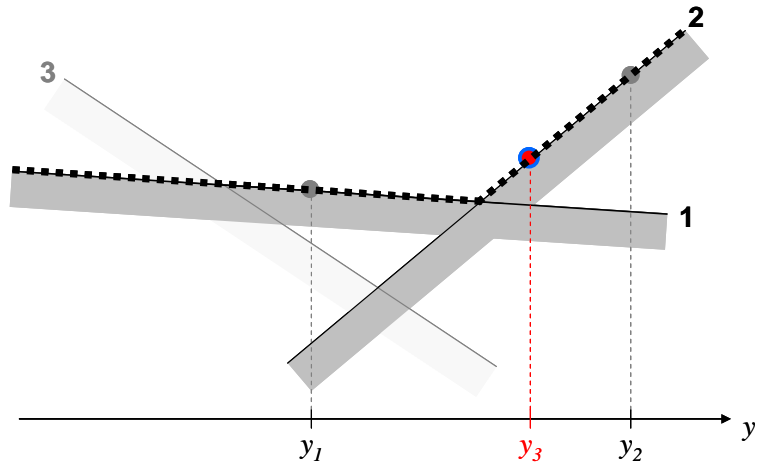

**Figure 17** *Step 3: Minimize the sub problem with $y_3$ as the value of the initial decision. No further constraint is added to the approximation of $Q(y, u_\omega)$.*

Here the active constraint at $y_3$ is one of the previously added constraints (constraint 2). That means that the objective value the optimizer thinks it will get from the sub-problems (blue circle in Figure 17) coincides with the actual value it does get (red circle). This will always be the case when a previously added constraint is the active constraint at the new optimal point. The reason for this is that our current approximation of $Q(y)$ (the heavy dotted line in figure 9) now correctly describes the sub-problem objective function at $y_3$. Since the optimizer chose this initial decision as optimal when it thought that the objective from the sub-problems would be at least at point C, when it turns out that this value actually is at point C the initial decision turns out to be optimal in light of the true sub-problem objective and so must be the overall optimal solution.

In general, if the approximation to the sub-problem objective function is a tight bound on the true sub-problem objective function at the current optimal initial decision then that initial decision is the true optimal initial decision.

Note that here we have found the optimal solution without adding constraint 3 to the master problem, so we have saved the need for a constraint when compared to the full formulation. We hope that this will be the case on a much larger scale when we solve real problems and that the amount of constraints we save will justify the extra effort solving and re-solving simplified problems.

### Unicut and Multicut

There are, in fact, two ways of adding optimality cuts to the master problem. The one described above is slightly simpler and is known as *unicut* because one cut is added to the master problem each iteration. The second method is known as *multicut*. This method was introduced by Birge and Louveaux (1969) and the idea is that instead of introducing a single additional cut on each iteration we introduce one for each of the sub problems.

This is based on the observation that $Q(y)$ is composed of a probability weighted sum of the sub-problem objective terms $Q(y, u_\omega)$. If instead of including $Q(y)$ as a single entity in the overall problem objective, we include it in its sum form we can then introduce a variable corresponding to the part of the objective function coming from each of the sub-problems. We can then place bounds on each of these variables in a similar way to that we used for bounding $Q(y)$.

So, we replace the optimization problem (34) with

$$\min_y c'y + \sum_\omega Q(y, u_\omega)$$
$$\text{s.t.} \quad Ay \leq b \quad (48)$$
$$y \geq 0.$$

We then replace each element of this sum with a variable $\theta_\omega$, which we bound using information from the sub-problems as before but without combining the bounds from each sub-problem into a single bound as in (44).

The cuts generated are, for each $\omega$,

$$\theta_\omega \geq p_\omega (\pi_\omega^i)'(h - Ty). \quad (49)$$

Or, letting $g_\omega = p_\omega (\pi_\omega^i)'h$ and $G_\omega = p_\omega (\pi_\omega^i)'T$,

$$\theta_\omega \geq g_\omega - G_\omega y \quad (50)$$

Instead of generating a single optimality cut in a given iteration as in unicut, we will generate one for each sub-problem. Birge and Louveaux (1988) show that multicut will converge more quickly than unicut for some problems. However, as the size of the master problem increases more quickly the storage requirements and solution time for each iteration are also increased.

### Algorithm

We can now present our final version of the Benders algorithm for stochastic programming problems, filling in the missing details from before.

Let *i* be the current iteration number, *r* be the number of feasibility cuts added and *s* be the number of optimality cuts added.

    0. Initialization

       a. Set $i = 0$, $r = 0$ and $s = 0$

1. Solve Master Problem

    a. Set $i = i + 1$

    b. If $i = 0$ solve the first stage problem only, i.e.
    $$\min_y c'y$$
    s.t.    $Ay \leq b$
    $$y \geq 0$$

    c. Otherwise, if $i > 0$, solve the current master problem:
    $$\min_y c'y + \theta$$
    s.t.    $Ay \leq b$
    $$y \geq 0$$
    $$D_j y \geq d_j \qquad j = 1,\ldots,r \qquad \text{(feasibility cuts)}$$
    $$\theta \geq g_j - G_j y \qquad j = 1,\ldots,s \qquad \text{(optimality cuts)}$$

    d. Let $y_i$ be the solution to the problem solved in b or c. Let $\theta_i$ be the value of $\theta$ calculated when solving this problem.

2. Feasibility Cuts

    a. Solve the feasibility variants of the sub-problems shown in (27) with the current solution of the master problem $y_i$.

    b. If any sub-problem is infeasible with this $y_i$, generate the feasibility cut in (29) / (30), increment $r$ and return to step 1.

3. Optimality Cuts

    a. All sub-problems are feasible at this point, since otherwise feasibility cuts would have been added above. For each of the sub-problems, solve its dual (40) and generate the optimality cut in (45) / (46) (for unicut) or the optimality cuts in (49) / (50) (for multicut).

4. Termination Test

    a. For each of the optimality cuts generated in step 3 above, test whether the current value of $\theta_i$ already satisfies the optimality cut conditions (for $i > 0$, as $\theta_0$ is not defined; if $i = 0$ do not terminate). If all optimality cut conditions are satisfied terminate, otherwise return to step 1 and re-solve.

In the worst case the feasibility and optimality cuts will recreate the original problem, but it is hoped that in practice the algorithm will terminate long before this happens.

There are some proven results about the convergence of solution methods based on Benders decomposition. A proof of convergence for bounded, feasible, two-stage problems is given by Van Slyke and Wets (1969) and Dempster and Merkovsky (1995) give a convergence proof for the case when $Q$ is piecewise and strictly convex. The theoretical worst-case convergence rate of the algorithm is geometric in the size of the problem, however in practice the algorithm appears to usually perform much better, though of course pathological examples can be constructed.

## Nested Benders

The development of Benders decomposition above was for two stage stochastic programming problems where we had an initial decision to make followed by a second-stage recourse decision at some later time, at which point we knew more information. In practice many stochastic programming problems consist of more than two decision stages. In such *multi-stage* problems an initial decision occurs followed by several stages of increasing information (e.g. with the passage of time) and further decision making. A good way to think of such problems is as a tree, with the initial decision at its root and each subsequent decision stage being a step further down the tree towards the leaves. The branching of the tree reflects the uncertainty in the development of the state of the world; on any particular branch you can arrive at any of the leaves attached to that branch. In discretized models the number of branches will be finite at each node, reflecting the discretization. As time goes by and the tree is descended, the number of possible leaves at which one can arrive is reduced, reflecting the increase of information (and consequent reduction in uncertainty) with the passage of time.

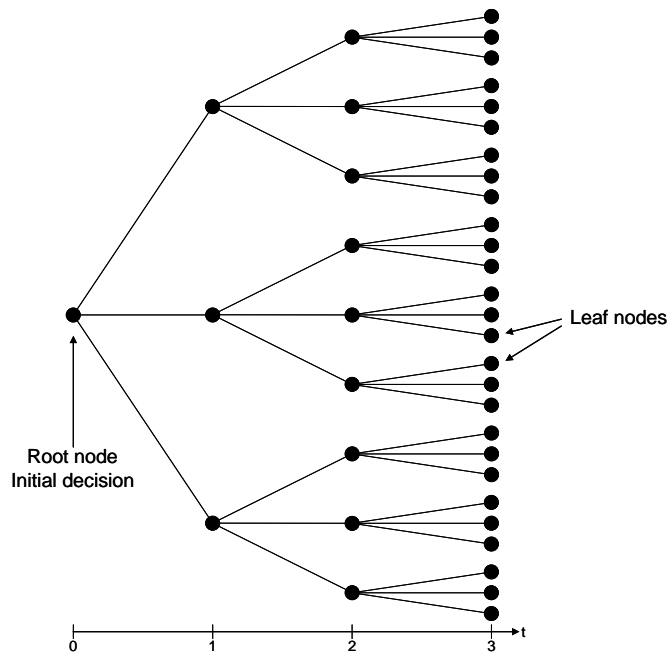

**Figure 18** *A 3.3.3 tree with four stages. Black circles represent decision points. Every decision point corresponds to a sub-problem.*

The *Nested Benders* algorithm gives us a way of applying the Benders decomposition method we have developed above for two stage problems to multi-stage problems. (It is sometimes also known as *L-shaped decomposition*, after the work of Van Slyke and Wets (1969)).

**Idea**

Nested Benders, as the name hints, takes the idea of Benders decomposition and applies it recursively over a tree structure, viewing the tree as set of nested two-stage problems.

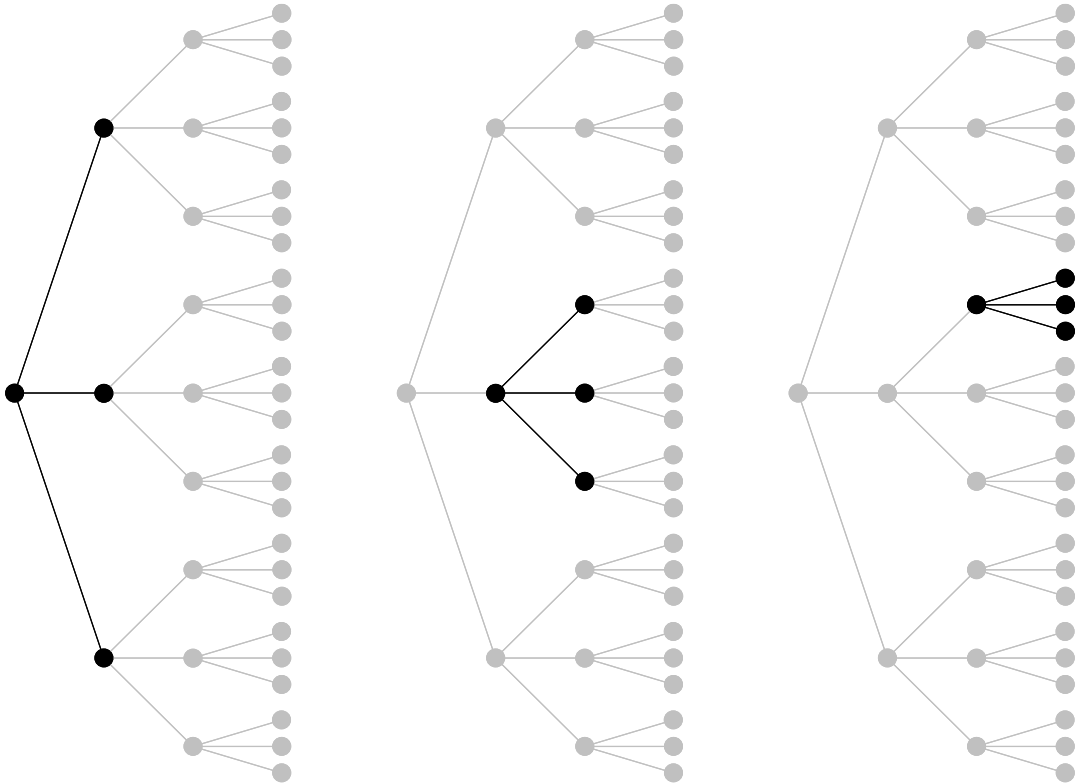

**Figure 19** *Viewing a multi-stage problem as a nested series of two-stage problems to which we can apply Benders decomposition. (a) left, the root node and successor nodes, (b) middle, the successor node and its children and (c) right, the third-level node and its leaf-node children.*

Looking at any pair of stages in isolation we have precisely the same problem as we dealt with in the preceding sections on Benders decomposition. The idea, then, is to apply Benders decomposition to each of these pairs of stages in turn, starting at the root node (Figure 19a). As we saw, solving decomposed problems relied on trying to solve the sub-problems with some chosen initial decision. At the root node, the sub-problems are those at each of the successor nodes at $t = 1$. So, we must try to solve each of these with the root-node decision already fixed. At this point, we can treat the initial decision as a constant and view the successor node with this particular fixed initial decision as the root node of its own problem (Figure 19b). This can continue down the tree until we get all the way to the leaves (Figure 19c).

In Benders decomposition we solve the sub-problems with the current initial decision and use these solutions to create feasibility and optimality cuts in the root-node master problem. In Nested Benders we do exactly the same at each level of the tree.

Using decomposition we can effectively isolate the problem at each stage. At each stage (other than the root), the connection to the previous stage is made through being passed a fixed decision from the parent stage and through the cuts sent back to the parent stage, which enforce feasibility or approximate the objective function. This means that we can solve parts of the problem at any stage we like as long as we have some value for the preceding decision.

We proceed as follows. First we solve the root node problem (with no constraints from subsequent stage sub-problems). We then take this optimal root decision and solve the sub-problems at the second stage (without any constraints coming from their subsequent stage sub-problems). These will generate feasibility or optimality cuts. So far, so like Benders. After solving the second stage problems we now have a choice. We can go back and solve the root node problem again using the cut or cuts that have been generated from the second stage problems or we can take the second stage optimal solutions that we have found (for the second stage problems that proved to be feasible) and move to the third stage, solving the sub-problems there. If we go to the third stage we will solve the sub-problems there (again without any constraints coming from their sub-problems) and generate some cuts for the second stage problems. We could then, if we wanted, move back up to the second stage and re-solve there or move on to the fourth stage. This process can (if desired) continue until we get to the leaves of the tree when it is no longer possible to move any further down and we have to start moving back up the tree. A nice way to think of the algorithm is that solutions are passed down the tree from parent nodes to children whilst cuts are passed back up the tree from children to parents.

The obvious question from the above description is when should we go up the tree and when should we go down it? This question has to be decided heuristically. That is to say, it has no correct answer, though it might well have wrong ones.

There is one obvious case. When all the sub-problems at a particular stage are infeasible it is impossible to descend the tree any further as we do not have solutions that we can take down the tree to the next stage. At this point we have to take the feasibility cuts generated and go back up the tree (though how far back up we go is again a heuristic).

One further thing that is helpful to us as we progress through the algorithm is that it is only necessary to re-solve sub-problems that have been modified since the last time we encountered them. This means we only need to solve (or re-solve) sub-problems that a) we have not encountered before, b) are being passed a new previous-stage decision from their parent node, or c) have had new cut information added to them from their children nodes.

We cannot descend further down any branches where we have encountered infeasibility (since in this case we have no proposed solution to pass to child nodes), or beyond the leaf nodes, and we cannot go up beyond the root node. These are, in fact, the only restrictions on our strategy for moving up and down the tree in the Nested Benders algorithm.

**Sequencing Heuristics**

There are three main heuristics that are used in deciding which way to move and when in the Nested Benders algorithm (Thompson, 1997). Such heuristics are known as *sequencing heuristics* or *sequencing protocols*.

- *Fast-forward*, which aims to descend the tree whenever possible. So, at any point where a sub-problem is feasible and we get a solution we will descend to the children nodes taking the solution with us until we either reach the leaves of the tree or infeasibility on all sub-problems under consideration.

- *Fast-back*, which aims to go back up the tree as far and as fast as it can wherever possible. So, whenever we add a cut from a child problem we will re-solve the parent and, if this generates a new cut, go to its parent and re-solve that, until we reach the root.

- *Fast-forward-fast-back*, which, as the name suggests, is a combination of the two methods above. In this case we have a switch that puts the algorithm in either descending or ascending mode. When in descending mode we follow the fast-forward scheme, going down the tree until we can go no further. When we can go no further we switch to ascending mode and then follow the fast-back scheme, ascending the tree until again we can go no further. We then switch back to descending mode.

This latter scheme has been found to be faster for a range of problems by Wittrock (1985) (and apparently Gassman (1987), though I have not seen his thesis) and has been used in practice. However, these are all heuristics; none have been shown to have theoretical advantages over the others and it should be possible to construct pathological examples that cause the performance of any one of these to be awful.

**Termination**

Termination of the Nested Benders algorithm is determined in the same way as for the original Benders algorithm. First, all the sub-problems at all levels of the tree have to be feasible. This means that all sub-problems have to be visited at least once and prove to be feasible with the solution proposed. If at any point the root node problem becomes infeasible then we must terminate as the entire problem is infeasible.

Assuming that the problem is feasible, a point will be reached where all the sub-problems are feasible with the current proposed solution. To terminate we require that all the sub-problems have reached an optimal solution. We saw in the section on Benders how to tell if, for a given problem, the solution is optimal. This was done by seeing whether the optimal value for the approximation of the sub-problem objective $\theta$ was greater or equal to the objective value coming from the sub-problems with the proposed solution $Q(y)$. (Remember that we are trying to minimize the objective, so if the chosen solution $y$ was thought optimal even when the value expected to be obtained from the sub-problem objectives was greater or equal to that which really was obtained then it is optimal – see the *Termination* sub-section from the previous section). In the case of Nested Benders this is an idea we can apply throughout the tree. If this test is satisfied for every sub-problem at all levels then the current proposed solution is optimal.

So, once we have a proposed solution *y* for a problem at a node we can see if it is optimal by taking the current value of $\theta$ (the sub-problem objective approximation) at the node and checking to see if it is greater or equal to the sub-problem objective function with the proposed parent node solution $Q(y)$. If this is the case at all nodes in the tree then the solution is optimal. This works because at each level of the tree we can be sure that the solution is optimal given the estimate of the objective functions from the sub-problems, which, if the test is passed, is at least as large as the true value.

### Algorithm

We will now present the Nested Benders algorithm.

0. Initialization
    a. Set current stage $t := 0$, the root stage
1. Solve the current problems at level $t$. It is only necessary to re-solve any problems that have not previously been solved or have had new cut information added to them.
    a. If any are infeasible place a feasibility cut in the parent problem at $t{-}1$. If $t = 0$ (root node) terminate as the problem is infeasible. Do not advance any further down infeasible branches. If all the feasible problems in this level have been solved to optimality (i.e. their children can add no further cut information with the current decision from the preceding stage) then ascend to parent, otherwise there is then a choice.
        i. Ascend ($t := t{-}1$) – return to the parent problem at level $t{-}1$ and re-solve (return to step 1).
        ii. Descend ($t := t{+}1$) – move down to the children of any feasible sub-problems and solve them (return to step 1). If all sub-problems are infeasible at level $t$ or we are at the leaf nodes ($t = T$) this option cannot be taken and we must ascend instead.
    b. There are no infeasible sub-problems at level $t$ (since otherwise step a. would have been followed). Place an optimality cut in the parent problem at $t{-}1$. If all the problems at this level have been solved before and will not be re-solved then descend (option ii. below), since no further cut information has been added to the parent. Otherwise there is then a choice.
        i. Ascend ($t := t{-}1$) – return to the parent problem at level $t{-}1$ and re-solve. This option cannot be taken if we are at the root level ($t = 0$). Instead move to step 2.

ii. Descend ($t := t+1$) – move down to the children of the sub-problems and solve them. If we are at the leaf nodes this option cannot be taken and we must ascend instead.

2. Termination test
   a. We are at the root node. For every sub-problem excluding the leaves (as they do not have children) check that the sub-problem objective approximation $\theta_X$ is greater or equal to the actual sub-problem objective term $Q(y_X)$, where here $X$ represents the particular sub-problem being examined.
   b. If $\theta_X \geq Q(y_X)$ for every sub-problem with children $X$ then terminate. We have reached optimality. If not then move down a level ($t := t+1$) and return to step 1.

### Algorithm Walkthrough

The operation of the Nested Benders algorithm is not (to me at least) completely obvious and so it might be helpful to see a simple example of it in operation. Consider applying the algorithm to the following tree.

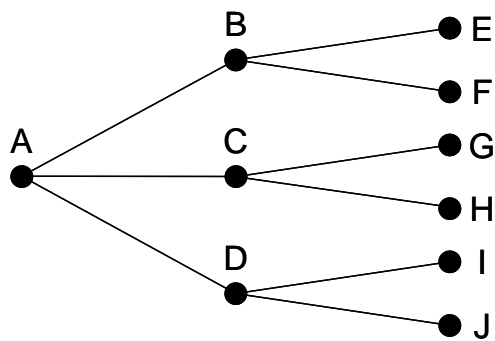

**Figure 20** *Tree for Nested Benders example walkthrough*

The algorithm will proceed as follows. Results, where needed, are made up and no particular sequencing heuristic is followed, though the possible choices of when to ascend or descend are noted.

- Solve root-node problem at A with no cuts
  - A feasible
- (At root node so can only move down tree)
- Solve B, C, D with no cuts
  - B, C feasible
  - D infeasible
- Place infeasibility cut in A from D
- Choice: Solve E, F and G, H or re-solve A

- - o   Fast-forward (FF) and Fast-forward-fast-back (FFFB) would solve E, F and G, H; Fast-back would solve A
- Solve E, F
  - o   Both feasible
- Place optimality cut from E, F into B
- Solve G, H
  - o   G feasible, H infeasible
- Place feasibility cut in C from H
- (Reached leaf nodes, can't descend further so go up)
- Re-solve B, C, D
  - o   Don't re-solve D as known to be infeasible with current initial decision
  - o   Do re-solve B, C as have new cut information in them
- Choice: Solve E, F and G, H or re-solve A
  - o   FF would solve EF, GH; FB or FFFB would solve A
  - o   For illustration choose to solve EF, GH.
- Solve E, F
  - o   Both feasible
- Place optimality cut from EF into B
  - o   Note that optimality cut placed in B is as previous one from EF: $\theta_B \geq Q_B(y_B)$
- Solve G, H
  - o   Both feasible
- Place optimality cut from GH into C
- (At leaves so must go up)
- Re-solve B, C, D
  - o   Know that D is infeasible with current initial decision, so don't re-solve
  - o   Know that B is optimal with current initial decision, so don't re-solve
  - o   Re-solve C as new optimality cut was added to it from children
- Choice: Solve G, H or re-solve A
  - o   FF would solve GH; FB or FFFB would solve A
  - o   Choose to solve GH
- Solve G, H

- o Both feasible
- Place optimality cut from GH into C
  - o Note that this is as previous one from GH: $\theta_C \geq Q_C(y_C)$
- (At leaves so must go up)
- Re-solve B, C, D
  - o D known to be infeasible with current initial decision
  - o B, C known to be optimal with current initial decision
  - o Don't re-solve anything at this stage
- (All feasible nodes in current stage solved to optimality, have to ascend to parent)
- Re-solve A
- (At root, must go down)
- Re-solve B, C, D with new root node decision from A
  - o All feasible
- Place optimality cut from B, C, D into A
- Choice: Solve EF, GH and IJ or re-solve A
  - o FF and FFFB would solve EF, GH and IJ; FB would re-solve A
  - o Choose to re-solve A
- Re-solve A
- (At root, must go down)
- Re-solve B, C, D with new value from A
  - o All feasible
- Add optimality cut to A
  - o Note optimality cut is as previous one added: $\theta_A \geq Q_A(y_A)$
- (Don't need to re-solve A as there are no new cuts added to it)
- Could have optimality!
  - o Check $\theta_X \geq Q_X(y_X)$ for all $X \in \{A, B, C, D\}$ (leaf nodes don't have a sub-problem objective approximation $\theta$).
  - o Find that it doesn't hold for D, so don't have optimality and must continue
- Go down to B, C, D and re-solve
  - o Don't need to re-solve any as all solved before
- (Can't go back up to A as have added no new cuts, must go down to EF, GH and IJ)

- Re-solve EF, GH, IJ
    - Don't need to re-solve EF, GH as have already solved them with current solutions at B and C respectively.
    - Re-solve I, J (have not solved before) – find them feasible
- Add optimality cut from I, J into D
- (At leaves, so must go back up)
- Re-solve B, C, D
    - B, C unchanged so do not need to be re-solved
    - D feasible, so all stage 1 problems feasible
- Add optimality cut from B, C, D to A
    - Note optimality cut is as one previously added: $\theta_A \geq Q_A(y_A)$
- (Don't need to re-solve A as no new cuts)
- Could have optimality!
    - Check $\theta_X \geq Q_X(y_X)$ for all $X \in \{A, B, C, D\}$
    - Holds for all $X$ – we have optimality!
- TERMINATE
    - Optimal decision given by current solution $y_X$ at each of the nodes

## Unicut and Multicut

The version of the algorithm described so far is a unicut variant. When we add an optimality cut to the parent, it is a single optimality cut coming from a combination of all of the sub-problems. We could replace this with the multicut version of optimality cuts (described in the earlier *Unicut and Multicut* section).

In the algorithm described above we continue to solve sub-problems even when an earlier one was discovered to be infeasible. It would also be possible to stop solving sub-problems as soon as infeasibility was discovered and then either proceed down a previously discovered feasible branch or return to the parent with a single feasibility cut. This could be seen as the strictest sort of unicut.

On the other hand, the most flexible sort of multicut could add an optimality cut or a feasibility cut to the parent from each sub-problem, depending on whether it was feasible or not, respectively.

## Aggregation

Aggregation is a technique where several stages of a stochastic programming problem are combined into a single stage for the purpose of solution efficiency.

We saw earlier that any stochastic programming problem has a deterministic equivalent form, which turns the problem from a multi-stage problem to a single stage linear programming problem. We can apply this same idea to multi-stage sub-sections of the tree, and collapse them down into a larger single-stage problem. If we do this to all the problems at a certain level of the tree we can reduce the overall number of stages in the tree. Obviously this is done at the cost of increasing the complexity of the sub-problems at the stages where aggregation has been applied. The following diagram illustrates the idea.

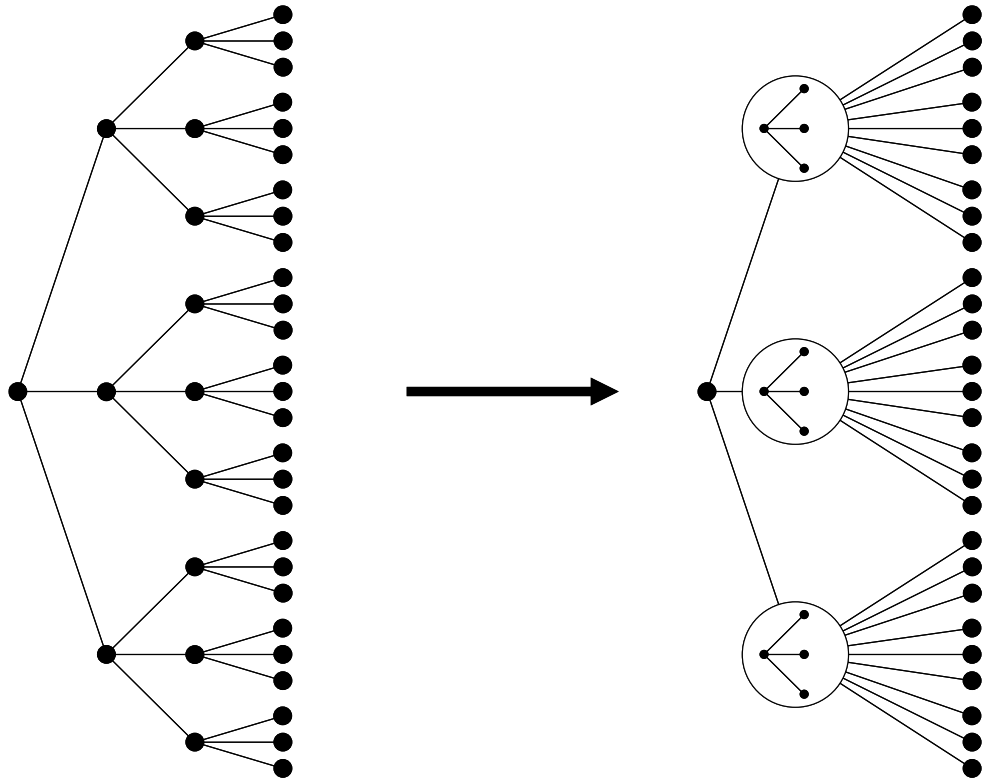

**Figure 21** *Aggregation of second and third stages to reduce tree from four to three stages. The large circles containing trees on the right represent the deterministic equivalent form of the stage 2-3 problems.*

Aggregation, along with our next topic, parallelization, is dealt with by Dempster and Thompson (1998), where they concluded that some aggregation could improve solution times when the nodal sub-problems were very small. Excessive aggregation resulted in increased solution times in their tests (presumably as the ability to exploit the problem structure is reduced). In that paper they suggest an optimal nodal problem size (of 100-300 rows and columns in the constraint matrix) and the reader is referred to that paper for further information.

**Parallelization**

It should be readily apparent that by breaking the task of solving a stochastic programming problem into smaller sub-problems that can each be solved somewhat autonomously, Benders and Nested Benders decomposition methods lend themselves in an obvious way to parallelization.

We will not cover this in great detail here and refer the reader to Dempster and Thompson (1998) where an algorithm for the parallelization of the Nested Benders algorithm is presented (along with references to other relevant work).

Their algorithm basically consists of having a master process controlling the despatch of problems to other processors and the current level at which we are working in the tree. At any given level, this master process farms out the sub-problems to the processors and receives feasibility or optimality information back from them. Once it has received information from all the sub-problems it can place an appropriate cut (or cuts in the case of multicut) in the parent problem and decide, according to its sequencing heuristic, whether to move up or down in the tree. This process can continue (and even the optimality test can be largely farmed out to the processors if required) until optimality is reached.

Dempster and Thompson (1998) give some results for the level of speed up found with different numbers of processors for different problems. They find that initially as the number of processes increases the speed-up is fairly close to linear, but the performance improvements decline as the number of processors becomes large due to the formation of bottlenecks.

## Alternative Methods

The Nested Benders algorithm is not the only algorithm available for the solution of large-scale stochastic programming problems. Whether or not it is the best probably depends on the problem being tackled. No attempt will be made here to describe these alternative methods (these are whole articles in themselves), but we will list them, if only to give some idea of what else to look up when investigating the solution of stochastic programming problems (and constrained optimization problems more generally). The list below is taken from the first chapter of Thompson (1997).

- Simplex method and variants (e.g. dual simplex)
- Lagrangian methods
- Interior point methods
    - Primal-dual log barrier methods
- Decomposition methods
    - Cutting plane decomposition methods (includes Nested Benders)
    - Augmented Lagrangian decomposition methods
- Sampling methods
- Stochastic quasigradient methods
- Derivative decomposition methods

In his summary of these methods, Thompson (1997) concludes that only simplex, Nested Benders and interior point methods need to be considered for the solution of stochastic programming problems on sequential machines.

## Conclusion

This article has described Benders decomposition, its application to two-stage stochastic programming problems and its extension to multi-stage problems in the form of the Nested Benders algorithm. It is hoped that this description will prove helpful to others attempting to understand this difficult but useful algorithm.

## Contact

If you have any comments or questions, please contact me at J@mesMurphy.co.uk, or via my website www.james-murphy.net. I will be happy to receive any feedback or corrections.